%% file: paper.tex
\documentclass[10pt,twocolumn,letterpaper]{article}

\usepackage{cvpr}
\usepackage{times}
\usepackage{amssymb}
\usepackage[fleqn]{amsmath}
\usepackage{url}
\usepackage{multirow}
\usepackage{amsthm} 
\usepackage{algorithm}
\usepackage{algorithmic}
\usepackage{diagbox} 
\usepackage{bbm}
\usepackage{bm}
\usepackage{diagbox}
\usepackage{graphicx}
\usepackage{subcaption}
\usepackage{epstopdf}
\usepackage{booktabs}
\usepackage{pdfpages}

\usepackage{courier}
\usepackage{color}
\definecolor{colorpink}{RGB}{251,53,155}
\definecolor{colorblue}{RGB}{0,148,200}
\definecolor{colorgreen}{RGB}{0,150,0}
\usepackage[pagebackref=true,breaklinks=true,colorlinks=true,urlcolor=cyan,linkcolor=red,citecolor=colorgreen]{hyperref}

\def\beq{\begin{eqnarray}\begin{aligned}}
\def\eeq{\end{aligned}\end{eqnarray}}
\def\noi{\noindent}
\def\nn{\nonumber}
\def\la{\langle}
\def\ra{\rangle}
\newtheorem{lemma}{Lemma}
\newtheorem{theorem}{Theorem}

\newcommand{\bbb}[1]{\boldsymbol{\mathbf{#1}}}

\def\Xi{\bar{\bbb{X}}}

\def\P{\mathcal{T}}

\def\O{\mathcal{O}}

\def\sign{\rm sign}

\def\ghs{\hspace{0.1cm}} 

\def\cone{\textbf}
\def\ctwo{\underline}
\def\cthree{\sout}

\newcommand{\algsize}{\fontsize{9.5}{12}\selectfont}

\usepackage{rotating}
\usepackage{array}

\def\imgheisubgraph{0.14\textheight}

\def\fourfigwid{0.235\textwidth}
\def\fivefigwid{0.195\textwidth}
\def\objimghei{0.1\textheight}
\usepackage{alphalph}

\DeclareCaptionFont{tiny}{\tiny}

\def\cone{\textcolor[rgb]{0.9961,0,0}}
\def\ctwo{\textcolor[rgb]{0,0.5820,0.9}}
\def\cthree{\textcolor[rgb]{0,0.7,0}}

\cvprfinalcopy 

\ifcvprfinal\fi
\def\horizontaldistance{\kern2pt}

\begin{document}

\title{A Matrix Splitting Method for Composite Function Minimization}

\author{Ganzhao Yuan\thanks{ Sun Yat-sen University (SYSU), King Abdullah University of Science \& Technology (KAUST). Email: yuanganzhao@gmail.com.},~~~Wei-Shi Zheng\thanks{Sun Yat-sen University (SYSU). Email: wszheng@ieee.org.},~~~Bernard Ghanem\thanks{King Abdullah University of Science \& Technology (KAUST). Email: bernard.ghanem@kaust.edu.sa.}}

\maketitle

\begin{abstract}
\input{abstract}
\end{abstract}

\input{introduction.tex}
\input{proposed.tex}

\input{extension.tex}

\input{experiment.tex}

\input{conclusion.tex}

{
\bibliographystyle{ieee}
\bibliography{egbib}
}

\end{document}

%% file: abstract.tex
Composite function minimization captures a wide spectrum of applications in both computer vision and machine learning. It includes bound constrained optimization and cardinality regularized optimization as special cases. This paper proposes and analyzes a new Matrix Splitting Method (MSM) for minimizing composite functions. It can be viewed as a generalization of the classical Gauss-Seidel method and the Successive Over-Relaxation method for solving linear systems in the literature. Incorporating a new Gaussian elimination procedure, the matrix splitting method achieves state-of-the-art performance. For convex problems, we establish the global convergence, convergence rate, and iteration complexity of MSM, while for non-convex problems, we prove its global convergence. Finally, we validate the performance of our matrix splitting method on two particular applications: nonnegative matrix factorization and cardinality regularized sparse coding. Extensive experiments show that our method outperforms existing composite function minimization techniques in term of both efficiency and efficacy.

%% file: introduction.tex
\section{Introduction}
In this paper, we focus on the following composite function minimization problem:
\beq \label{eq:main}
\min_{\bbb{x}}~f(\bbb{x}) \triangleq q(\bbb{x}) + h(\bbb{x});~q(\bbb{x}) = \tfrac{1}{2}\bbb{x}^T\bbb{A}\bbb{x} + \bbb{x}^T\bbb{b}
\eeq
\noi where $\bbb{x} \in\mathbb{R}^n,~\bbb{b} \in \mathbb{R}^n$, $\bbb{A}\in \mathbb{R}^{n\times n}$ is a positive semidefinite matrix, $h(\bbb{x})$ is a piecewise separable function (\ie $h(\bbb{x})=\sum_{i=1}^n h(\bbb{x}_i)$) but not necessarily convex. Typical examples of $h(\bbb{x})$ include the bound constrained function and the $\ell_0$ and $\ell_1$ norm functions.

The optimization in (\ref{eq:main}) is flexible enough to model a variety of applications of interest in both computer vision and machine learning, including compressive sensing \cite{Donoho06}, nonnegative matrix factorization \cite{lee1999learning,lin2007projected,guan2012nenmf}, sparse coding \cite{lee2006efficient,aharon2006k,BaoJQS16,Quan2016CVPR}, support vector machine \cite{hsieh2008dual}, logistic regression \cite{yu2011dual}, subspace clustering \cite{elhamifar2013sparse}, to name a few. Although we only focus on the quadratic function $q(\cdot)$, our method can be extended to handle general composite functions as well, by considering a typical Newton approximation of the objective \cite{TsengY09,yuan2014newton}.

The most popular method for solving (\ref{eq:main}) is perhaps the proximal gradient method \cite{nesterov2013introductory,beck2009fast}. It considers a fixed-point proximal iterative procedure $\bbb{x}^{k+1}= \text{prox}_{\gamma h}(\bbb{x}^k - \gamma \triangledown q(\bbb{x}^k))$ based on the current gradient $\triangledown q(\bbb{x}^k)$. Here the proximal operator $\text{prox}_{\tilde{h}}(\bbb{a}) = \arg \min_{\bbb{x}} ~\frac{1}{2}\|\bbb{x}-\bbb{a}\|_2^2 + \tilde{h}(\bbb{x})$ can often be evaluated analytically, $\gamma={1}/{L}$ is the step size with $L$ being the local (or global) Lipschitz constant. It is guaranteed to decrease the objective at a rate of $\mathcal{O}({L}/{k})$, where $k$ is the iteration number. The accelerated proximal gradient method can further boost the rate to $\mathcal{O}({L}/{k^2})$. Tighter estimates of the local Lipschitz constant leads to better convergence rate, but it scarifies additional computation overhead to compute $L$. Our method is also a fixed-point iterative method, but it does not rely on a sparse eigenvalue solver or line search backtracking to compute such a Lipschitz constant, and it can exploit the specified structure of the quadratic Hessian matrix $\mathbf{A}$.

The proposed method is essentially a generalization of the classical Gauss-Seidel (GS) method and Successive Over-Relaxation (SOR) method \cite{demmel1997applied,saad2003iterative}. In numerical linear algebra, the Gauss-Seidel method, also known as the successive displacement method, is a fast iterative method for solving a linear system of equations. It works by solving a sequence of triangular matrix equations. The method of SOR is a variant of the GS method and it often leads to faster convergence. Similar iterative methods for solving linear systems include the Jacobi method and symmetric SOR. Our proposed method can solve versatile composite function minimization problems, while inheriting the efficiency of modern linear algebra techniques.


\begin{figure} [!t]
\centering
      \begin{subfigure}{\fourfigwid}\includegraphics[width=\textwidth,height=\imgheisubgraph]{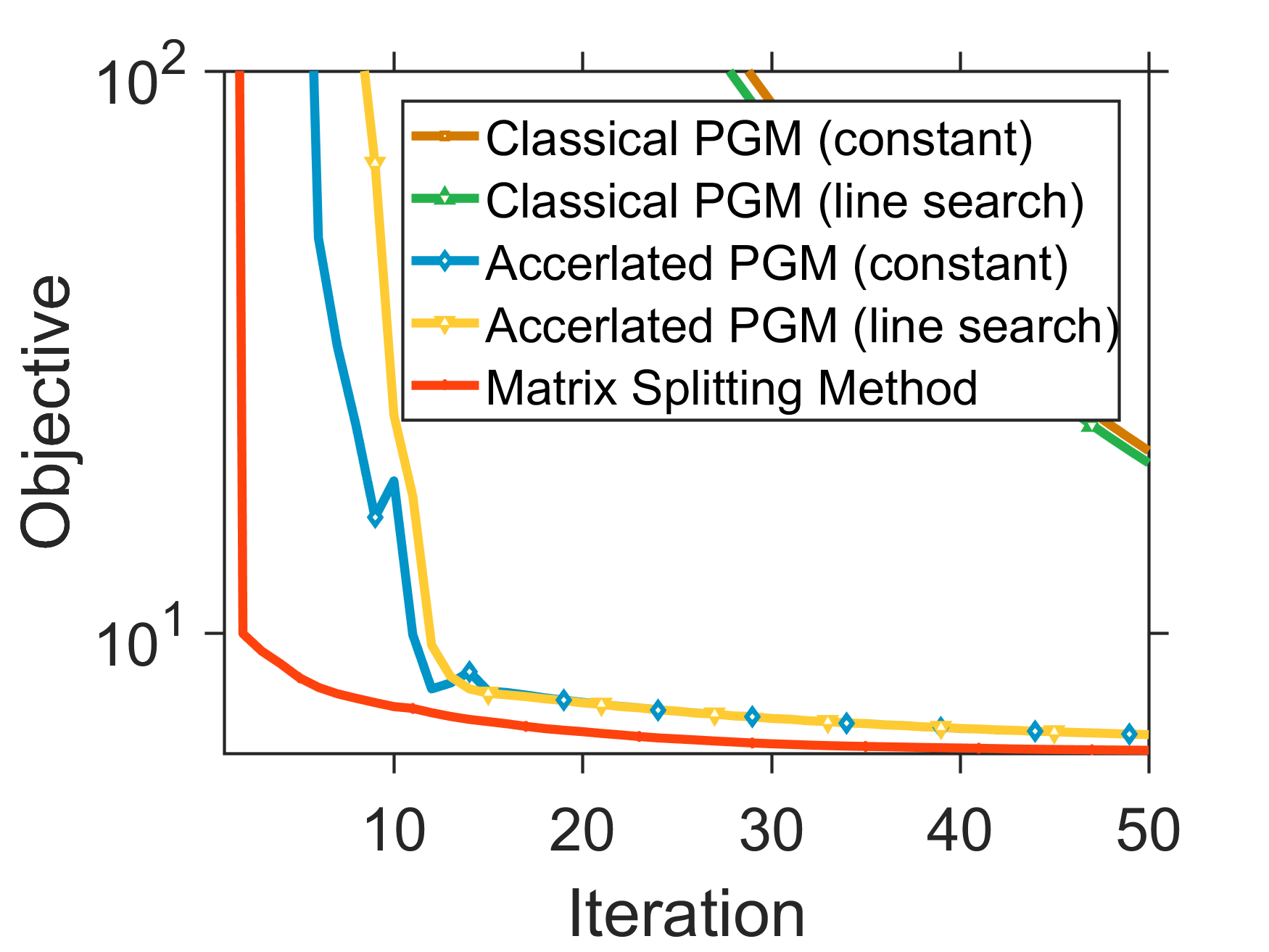}\vspace{-6pt} \caption{\footnotesize $\underset{\bbb{x}\geq \mathbf{0}}{\min}~~\tfrac{1}{2}\|\bbb{Cx}-\bbb{d}\|_2^2$}\end{subfigure}\ghs
      \begin{subfigure}{\fourfigwid}\includegraphics[width=\textwidth,height=\imgheisubgraph]{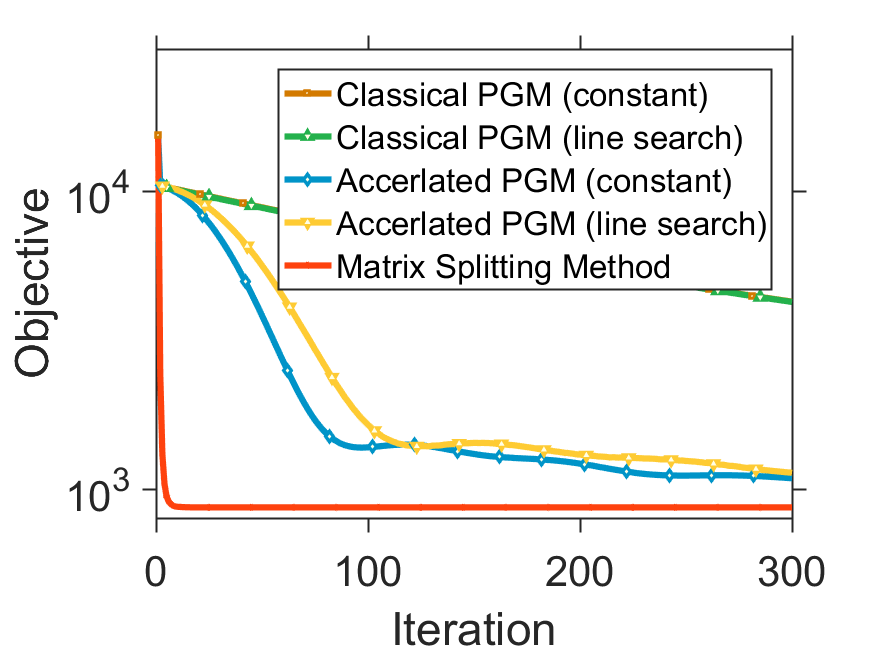}\vspace{-6pt} \caption{\footnotesize $\underset{\bbb{x}}{\min}~~\tfrac{1}{2}\|\bbb{Cx}-\bbb{d}\|_2^2+\|\bbb{x}\|_0$}\end{subfigure}
\vspace{-8pt}
\caption{Convergence behavior for solving (a) convex non-negative least squares and (b) nonconvex $\ell_0$ norm sparse least squares. We generate $\bbb{C}\in\mathbb{R}^{200\times 1000}$ and $\bbb{d}\in\mathbb{R}^{200}$ from a (0-1) uniform distribution. All methods share the same initial point. Our matrix splitting method significantly outperforms existing popular proximal gradient methods \cite{beck2009fast,nesterov2013introductory} such as classical Proximal Gradient Method (PGM) with constant step size, classical PGM with line search, accelerated PGM with constant step size and accelerated PGM with line search. Note that all the methods have the same computational complexity for one iteration.}
\label{fig:intro}
\vspace{-10pt}
\end{figure}

Our method is closely related to coordinate gradient descent and its variants such as randomized coordinate descent \cite{hsieh2008dual}, cyclic coordinate descent, block coordinate descent \cite{nesterov2012efficiency}, accelerated randomized coordinate descent \cite{lu2013randomized} and others \cite{richtarik2014iteration,liu2015asynchronous,sun2015improved,zeng2015gauss}. However, all these work are based on gradient-descent type iterations and a constant Lipschitz step size. They work by solving a first-order majorization/surrogate function via closed form updates. Their algorithm design and convergence result cannot be applied here. In contrast, our proposed method does not rely on computing the gradient descent direction or the Lipschicz constant step size, yet it adopts a triangle matrix factorization strategy, where the triangle subproblem can be solved by an alternating cyclic coordinate strategy.

\vspace{3pt}\noindent \textbf{Contributions.} \textbf{(a)} We propose a new Matrix Splitting Method (MSM) for composite function minimization. \textbf{(b)} For convex problems, we establish the global convergence, convergence rate, and iteration complexity of MSM, while for non-convex problems, we prove its global convergence to a local optimum. \textbf{(c)} Our experiments on nonnegative matrix factorization and sparse coding show that MSM outperforms state-of-the-art approaches. Before proceeding, we present a running example in Figure \ref{fig:intro} to show the performance of our proposed method, as compared with existing ones.

\vspace{3pt}\noindent \bbb{Notation.} We use lowercase and uppercase boldfaced letters to denote real vectors and matrices respectively. The Euclidean inner product between $\bbb{x}$ and $\bbb{y}$ is denoted by $\la \bbb{x},\bbb{y}\ra$ or $\bbb{x}^T\bbb{y}$. We denote $\|\bbb{x}\|=\|\bbb{x}\|_2=\sqrt{\la \bbb{x},\bbb{x} \ra}$, and $\|\bbb{C}\|$ as the spectral norm (\ie the largest singular value) of $\bbb{C}$. We denote the $i^{\text{th}}$ element of vector $\bbb{x}$ as $\bbb{x}_i$ and the $(i,j)^{\text{th}}$ element of matrix $\mathbf{C}$ as $\mathbf{C}_{i,j}$. $diag(\bbb{D}) \in \mathbb{R}^{n}$ is a column vector formed from the main diagonal of $\bbb{D}\in\mathbb{R}^{n\times n}$. $\bbb{C}\succeq0$ indicates that the matrix $\bbb{C}\in\mathbb{R}^{n\times n}$ is positive semidefinite (not necessarily symmetric)\footnote{$\bbb{C}\succeq0 \Leftrightarrow \forall \bbb{x},~\bbb{x}^T\bbb{Cx}\geq 0 \Leftrightarrow \forall \bbb{x},~\frac{1}{2}\bbb{x}^T(\bbb{C}+\bbb{C}^T)\bbb{x}\geq 0$}. Finally, we denote $\bbb{D}$ as a diagonal matrix of $\bbb{A}$ and $\bbb{L}$ as a strictly lower triangle matrix of $\bbb{A}$ \footnote{For example, when $n=3$, $\bbb{D}$ and $\bbb{L}$ take the following form: \\$\textstyle \renewcommand\arraystretch{1}
\setlength{\arraycolsep}{3pt}
\bbb{D}=\begin{bmatrix}
  \bbb{A}_{1,1} & 0 & 0    \\
  0 & \bbb{A}_{2,2} & 0  \\
  0 & 0 & \bbb{A}_{3,3}   \\
\end{bmatrix},~\bbb{L}=\begin{bmatrix}
  0 & 0 & 0    \\
  \bbb{A}_{2,1} & 0 & 0  \\
  \bbb{A}_{3,1} & \bbb{A}_{3,2} &0    \\
\end{bmatrix}$}. Thus, we have $\bbb{A} = \bbb{L}+\bbb{D}+\bbb{L}^T$.

%% file: proposed.tex
\section{Proposed Matrix Splitting Method}

In this section, we present our proposed matrix splitting method for solving (\ref{eq:main}). Throughout this section, we assume that $h(\bbb{x})$ is convex and postpone the discussion for nonconvex $h(\bbb{x})$ to Section \ref{sect:extension}.

Our solution algorithm is derived from a fixed-point iterative method based on the first-order optimal condition of (\ref{alg:main}). It is not hard to validate that a solution $\bbb{x}$ is the optimal solution of (\ref{alg:main}) if and only if $\bbb{x}$ satisfies the following nonlinear equation (``$\triangleq$'' means define):
\beq \label{eq:opt:cond}
\textstyle \bbb{0} \in \partial f(\bbb{x})\triangleq \bbb{Ax} + \bbb{b} + \partial h(\bbb{x})
\eeq
\noi where $\partial h(\bbb{x})$ is the sub-gradient of $h(\cdot)$ in $\bbb{x}$. In numerical analysis, a point $\bbb{x}$ is called a fixed point if it satisfies the equation $\bbb{x} \in \P(\bbb{x})$, for some operator $\P(\cdot)$. Converting the transcendental equation $\bbb{0}\in\partial f(\bbb{x})$ algebraically into the form $\bbb{x} \in \P(\bbb{x})$, we obtain the following iterative scheme with recursive relation:
\beq \label{eq:iterative}
\textstyle \bbb{x}^{k+1}\in \P(\bbb{x}^k),~k = 0, 1, 2,...
\eeq
We now discuss how to adapt our algorithm into the iterative scheme in (\ref{eq:iterative}). First, we split the matrix $\bbb{A}$ in (\ref{eq:opt:cond}) using the following strategy:
\beq\label{eq:matrix}
\textstyle  \bbb{A} = \underbrace{\bbb{L}+\tfrac{1}{\omega}(\bbb{D}+\theta \bbb{I})}_{\bbb{B}} + \underbrace{ \bbb{L}^T + \tfrac{1}{\omega}((\omega-1)\bbb{D} - \theta \bbb{I})}_{\bbb{C}}
\eeq
\noi Here, $\omega \in (0,2)$ is a relaxation parameter and $\theta\in[0,\infty)$ is a parameter for strong convexity that enforces $diag(\bbb{B})>\bbb{0}$. These parameters are specified by the user beforehand. Using these notations, we obtain the following optimality condition which is equivalent to (\ref{eq:opt:cond}):
\beq \label{eq:kkt:BC}
\textstyle \textstyle-\bbb{Cx} -  \bbb{b} \in  (\bbb{B}+\partial h ) (\bbb{x})   \nn
\eeq
\noi Then, we have the following equivalent fixed-point equation:
\beq\label{eq:P}
\textstyle \bbb{x} \in \P(\bbb{x}) \triangleq  (\bbb{B}+\partial h )^{-1}(-\bbb{Cx} -  \bbb{b})
\eeq
\noi Here, we name $\P$ the triangle proximal operator, which is novel in this paper\footnote{This is in contrast with Moreau's proximal operator \cite{parikh2014proximal}: $\text{prox}_{h}(\bbb{a}) = \arg\min_{\bbb{x}} ~\frac{1}{2}\|\bbb{x}-\bbb{a}\|_2^2 + h(\bbb{x})=(\bbb{I}+\partial h)^{-1}(\bbb{a})$, where the mapping $(\bbb{I}+\partial h)^{-1}$ is called the resolvent of the subdifferential operator $\partial h$.}. Due to the triangle property of the matrix $\bbb{B}$ and the element-wise separable structure of $h(\cdot)$, the triangle proximal operator $\P(\bbb{x})$ in (\ref{eq:P}) can be computed exactly and analytically, by a generalized Gaussian elimination procedure (discussed later in Section \ref{sect:alg:subprob}). Our matrix splitting method iteratively applies $\bbb{x}^{k+1} \Leftarrow\P(\bbb{x}^k)$ until convergence. We summarize our algorithm in Algorithm \ref{alg:main}.



In what follows, we show how to compute $\P(\bbb{x})$ in (\ref{eq:P}) in Section \ref{sect:alg:subprob}, and then we study the convergence properties of Algorithm \ref{alg:main} in Section \ref{sect:alg:convergence}.




\subsection{Computing the Triangle Proximal Operator}\label{sect:alg:subprob}

We now present how to compute the triangle proximal operator in (\ref{eq:P}), which is based on a new generalized Gaussian elimination procedure. Notice that (\ref{eq:P}) seeks a solution $\bbb{z}^*\triangleq \P(\bbb{x}^k)$ that satisfies the following nonlinear system:
\beq \label{eq:subproblem}
\textstyle \bbb{0} \in \bbb{B}\bbb{z}^*   +  \bbb{u} +  \partial h(\bbb{z}^*)   ,~\text{where}~\bbb{u}=\bbb{b} + \bbb{C}\bbb{x}^k
\eeq
\noi By taking advantage of the triangular form of $\bbb{B}$ and the element-wise structure of $h(\cdot)$, the elements of $\bbb{z}^*$ can be computed sequentially using forward substitution. Specifically, the above equation can be written as a system of nonlinear equations:

{
\tiny
\footnotesize
\renewcommand\arraystretch{1}
\setlength{\arraycolsep}{1.5pt}
\beq \label{eq:sub:nonlinear}
\bbb{0} \in \begin{bmatrix}
\bbb{B}_{1,1}& 0 & 0 & 0 &0  \\
\bbb{B}_{2,1}& \bbb{B}_{2,2} & 0 & 0 &0 \\
\vdots & \vdots  & \ddots &    0 &   0 \\
\bbb{B}_{n-1,1} &\bbb{B}_{n-1,2} & \cdots & \bbb{B}_{n-1,n-1}  & 0\\
\bbb{B}_{n,1}& \bbb{B}_{n,2} & \cdots & \bbb{B}_{n,n-1} &\bbb{B}_{n,n}  \\
\end{bmatrix} \begin{bmatrix}
  \bbb{z}^*_{1} \\
  \bbb{z}^*_{2}\\
  \vdots \\
  \bbb{z}^*_{n-1}   \\
  \bbb{z}^*_{n}   \\
\end{bmatrix} + \bbb{u} + \partial h(\bbb{z}^*) \nn
\eeq
}
\noi If $\bbb{z}^*$ satisfies the equations above, it must solve the following one-dimensional subproblems:
\beq
0 \in \bbb{B}_{j,j} \bbb{z}^*_j + \bbb{w}_j +  \partial  h(\bbb{z}^*_j),~\forall j=1,2,~...~,n,\nn\\
\textstyle \bbb{w}_j=\bbb{u}_j +  \sum_{i=1}^{j-1} \bbb{B}_{j,i}\bbb{z}^*_{i}~~~~~~~~~~~~~~~~~\nn
\eeq
\noi This is equivalent to solving the following one-dimensional problem for all $j=1,2,...,n$:
\beq\label{eq:1d:subp}
\textstyle \bbb{z}^*_j= t^* \triangleq \underset{t}{\arg\min}~~\tfrac{1}{2}\bbb{B}_{j,j} t^2 + \bbb{w}_j t +  h(t)
\eeq
\noi Note that the computation of $\bbb{z}^{*}$ uses only the elements of $\bbb{z}^{*}$ that have already been computed and a successive displacement strategy is applied to find $\bbb{z}^{*}$.

We remark that the one-dimensional subproblem in (\ref{eq:1d:subp}) often admits a closed form solution for many problems of interest. For example, when $h(t)=I_{[lb,ub]}(t)$ with $I(\cdot)$ denoting an indicator function on the box constraint $lb\leq \bbb{x}\leq ub$, the optimal solution can be computed as: $t^* = \min(ub,\max(lb,-\bbb{w}_j/\bbb{B}_{j,j}))$; when $h(t)=\lambda |t|$ (\eg in the case of the $\ell_1$ norm), the optimal solution can be computed as: $t^* =   - \max\left(0,| \bbb{w}_j/\bbb{B}_{j,j}|-\lambda/\bbb{B}_{j,j}\right) \cdot \sign\left(\bbb{w}_j/\bbb{B}_{j,j}\right) $.

Our generalized Gaussian elimination procedure for computing $\P(\bbb{x}^k)$ is summarized in Algorithm \ref{alg:sub}. Note that its computational complexity is $\O(n^2)$, which is the same as computing a matrix-vector product.

\begin{algorithm} [!t]
\algsize
\caption{\label{alg:main} {\bbb{MSM}: A Matrix Splitting Method for Solving the Composite Function Minimization Problem in (\ref{eq:main})}}
\begin{algorithmic}[1]
\STATE Choose $\omega\in(0,2),~\theta\in[0,\infty)$.~Initialize $\bbb{x}^0$, ~$k=0$.\\
\STATE \text{while not converge}\\
\STATE~~~$\bbb{x}^{k+1} = \P(\bbb{x}^k)$ (Solve (\ref{eq:subproblem}) by Algorithm \ref{alg:sub})
\STATE~~~$k = {k+1}$
\STATE \text{end while}\\
\STATE Output $\bbb{x}^{k+1}$\\
\end{algorithmic}
\end{algorithm}
\begin{algorithm} [!t]
\algsize
\caption{\label{alg:sub} {A Generalized Gaussian Elimination Procedure for Computing the Triangle Proximal Operator $\P(\bbb{x}^k)$.}}
\begin{algorithmic}[1]
\STATE Input $\bbb{x}^k$\\
\STATE Initialization: compute $\bbb{u}= \bbb{b} + \bbb{C}\bbb{x}^k$\\
\STATE $\bbb{x}_1=\arg\min_{t}\frac{1}{2}\bbb{B}_{1,1}t^2 + (\bbb{u}_1) t +  h(t)$\\
\STATE $\bbb{x}_2=\arg\min_{t}\frac{1}{2}\bbb{B}_{2,2}t^2 + (\bbb{u}_2   + \bbb{B}_{2,1}\bbb{x}_{1}) t +  h(t)$\\
\STATE $\bbb{x}_3=\arg\min_{t}\frac{1}{2}\bbb{B}_{3,3}t^2 + (\bbb{u}_3+ \bbb{B}_{3,1}\bbb{x}_{1} + \bbb{B}_{3,2}\bbb{x}_{2}) t +  h(t)$\\
\STATE ...\\
\STATE $\bbb{x}_{n}=\arg\min_{t}\frac{1}{2}\bbb{B}_{n,n}t^2 + (\bbb{u}_n +  \sum_{i=1}^{n-1} \bbb{B}_{n,i}\bbb{x}_{i}) t +  h(t)$\\
\STATE Collect $(\bbb{x}_1,\bbb{x}_2,\bbb{x}_3,...,\bbb{x}_n)^T$ as $\bbb{x}^{k+1}$ and Output $\bbb{x}^{k+1}$\\
\end{algorithmic}
\end{algorithm}
%


\subsection{Convergence Analysis}\label{sect:alg:convergence}
In what follows, we present our convergence analysis for Algorithm \ref{alg:main}. We let $\bbb{x}^*$ be the optimal solution of (\ref{eq:main}). For notation simplicity, we denote:
\beq
&~~~~~~~~~\bbb{r}^k\triangleq \bbb{x}^{k} - \bbb{x}^*,~\bbb{d}^k \triangleq \bbb{x}^{k+1}-\bbb{x}^k \\
&u^k \triangleq f(\bbb{x}^{k})-f(\bbb{x}^*),~f^k \triangleq f(\bbb{x}^k),~f^* \triangleq f(\bbb{x}^*)
\eeq

\noindent The following lemma characterizes the optimality of $\P(\bbb{y})$.

\begin{lemma} \label{lemma:opt:ineq}
For all $\bbb{x},\bbb{y}\in\mathbb{R}^n$, it holds that:
\beq
\bbb{v}\in \partial h(\P(\bbb{x})),~\|\bbb{A}\P(\bbb{x}) + \bbb{b} +  \bbb{v}\| \leq  \|\bbb{C}\| \|\bbb{x}- \P(\bbb{x})\|   \label{eq:opt:bound0}
\eeq
\vspace{-18pt}
\beq\label{eq:opt:bound}
&  f(\P(\bbb{y}))-  f(\bbb{x}) \leq   \la \P(\bbb{y})-\bbb{x},\bbb{C}(\P(\bbb{y})-\bbb{y})  \\
&~~~~~~~~~~ - \tfrac{1}{2}(\bbb{x}-\P(\bbb{y}))^T\bbb{A}(\bbb{x}-\P(\bbb{y}))  \ra
\eeq
\begin{proof}
(i) We now prove (\ref{eq:opt:bound0}). By the optimality of $\P(\bbb{x})$, we have: $\forall \bbb{v}\in \partial h(\P(\bbb{x})),~0 = \bbb{B}\P(\bbb{x}) + \bbb{C}\P(\bbb{x}) +  \bbb{v} + \bbb{b} + \bbb{C} \bbb{x} - \bbb{C}\P(\bbb{x})$. Therefore, we obtain: $\bbb{A}\P(\bbb{x})  +  \bbb{b}+ \bbb{v} =  \bbb{C}(\P(\bbb{x})-\bbb{x})$. Applying a norm inequality, we have (\ref{eq:opt:bound0}).

(ii) We now prove (\ref{eq:opt:bound}). For simplicity, we denote $\bbb{z}^* \triangleq \P(\bbb{y})$. Thus, we obtain: $\bbb{0} \in ~\bbb{B} \bbb{z}^* +  \bbb{b} + \bbb{C}\bbb{y} + \partial h(\bbb{z}^*)\Rightarrow 0 \in ~\la\bbb{x}-\bbb{z}^*, \bbb{B} \bbb{z}^* + \bbb{b} + \bbb{C}\bbb{y} + \partial h(\bbb{z}^*)  \ra,~\forall \bbb{x}$.
\noi Since $h(\cdot)$ is convex, we have:
\beq \label{eq:fail}
\textstyle \la\bbb{x}-\bbb{z}^*,\partial h(\bbb{z}^*)  \ra \leq  h(\bbb{x}) - h(\bbb{z}^*)
\eeq
Then we have this inequality: $\forall \bbb{x}:~h(\bbb{x}) -  h(\bbb{z}^*) + \la\bbb{x}-\bbb{z}^*, \bbb{B} \bbb{z}^* + \bbb{b} +\bbb{C} \bbb{y} \ra  \geq 0$. We naturally derive the following results: $f(\bbb{z}^*)- f(\bbb{x}) = h(\bbb{z}^*)  -  h(\bbb{x}) +  q(\bbb{z}^*) -  q(\bbb{x})\leq \la \bbb{x}-\bbb{z}^*, \bbb{B}\bbb{z}^* +  \bbb{b} + \bbb{C}\bbb{y}  \ra  +  q(\bbb{z}^*) -  q(\bbb{x}) =\la \bbb{x}-\bbb{z}^*, \bbb{B}\bbb{z}^* + \bbb{C}\bbb{y}  \ra + \tfrac{1}{2}{\bbb{z}^*}^T\bbb{A}\bbb{z}^* -\tfrac{1}{2}\bbb{x}^T\bbb{A}\bbb{x}=\la \bbb{x}-\bbb{z}^*, (\bbb{B}-\bbb{A})\bbb{z}^* + \bbb{C}\bbb{y}  \ra -\tfrac{1}{2} (\bbb{x}-\bbb{z}^*)^T \bbb{A} (\bbb{x}-\bbb{z}^*)=\la \bbb{x}-\bbb{z}^*, \bbb{C}(\bbb{y}-\bbb{z}^*) \ra -\tfrac{1}{2} (\bbb{x}-\bbb{z}^*)^T \bbb{A} (\bbb{x}-\bbb{z}^*)$.

\end{proof}
\end{lemma}

\begin{theorem} \label{theorem:1}
(Proof of Global Convergence) We define $\delta \triangleq \tfrac{2\theta}{\omega} + \tfrac{2-\omega}{\omega}  \min(diag(\bbb{D}))$. Assume that $\omega$ and $\theta$ are chosen such that $\delta\in(0,\infty)$, Algorithm \ref{alg:main} is globally convergent.

\begin{proof}
(i) First, the following results hold for all $\bbb{z}\in\mathbb{R}^n$:
\beq \label{eq:upperbound:0}
 \bbb{z}^T(\bbb{A}-2\bbb{C})\bbb{z}
 =~&\textstyle \textstyle \bbb{z}^T(\bbb{L} - \bbb{L}^T  + \tfrac{2-\omega}{\omega}\bbb{D} + \tfrac{2\theta}{\omega}  \bbb{I})\bbb{z} \\
=~&\textstyle \bbb{z}^T( \tfrac{2\theta}{\omega} \bbb{I} + \tfrac{2-\omega}{\omega}\bbb{D} )\bbb{z} \geq  \delta \|\bbb{z}\|_2^2
\eeq
\noi where we have used the definition of $\bbb{A}$ and $\bbb{C}$, and the fact that $\bbb{z}^T \bbb{L}\bbb{z} = \bbb{z}^T\bbb{L}^T\bbb{z},~\forall \bbb{z}$.

We invoke (\ref{eq:opt:bound}) in Lemma \ref{lemma:opt:ineq} with $\bbb{x}=\bbb{x}^k,~\bbb{y}=\bbb{x}^k$ and combine the inequality in (\ref{eq:upperbound:0}) to obtain:
\beq \label{eq:descent}
 \textstyle f^{k+1} - f^k & \textstyle \leq  -\frac{1}{2} \la \bbb{d}^{k},(\bbb{A}-2\bbb{C})\bbb{d}^{k}\ra \leq \textstyle -\frac{\delta}{2} \|\bbb{d}^k\|_2^2
\eeq

\noi (ii) Second, we invoke (\ref{eq:opt:bound0}) in Lemma \ref{lemma:opt:ineq} with $\bbb{x}=\bbb{x}^k$ and obtain: $\bbb{v}\in \partial h(\bbb{x}^{k+1}),~\tfrac{1}{ \|\bbb{C}\|} \|\bbb{A}\bbb{x}^{k+1} + \bbb{b} +  \bbb{v}\| \leq  \|\bbb{x}^k-\bbb{x}^{k+1}\|$. Combining with (\ref{eq:descent}), we have: $\frac{\delta }{2\|\bbb{C}\|} \cdot\|\partial f(\bbb{x}^{k+1})\|_2^2 \leq f^k - f^{k+1}$, where $\partial f(\bbb{x}^k)$ is defined in (\ref{eq:opt:cond}). Summing this inequality over $i=0,...,k-1$, we have: $\frac{\delta }{2\|\bbb{C}\|} \cdot \sum_{i=0}^{k-1} \|\partial f(\bbb{x}^i)\|_2^2 \leq \textstyle f^0  - f^{k}\leq \textstyle  f^0  - f^* \nn$, where we use $f^*\leq f^{k}$. As $k\rightarrow \infty$, we have $\partial f(\bbb{x}^k)\rightarrow \bbb{0}$, which implies the convergence of the algorithm.

Note that guaranteeing $\delta\in(0,\infty)$ can be achieved by simply choosing $\omega\in(0,2)$ and setting $\theta$ to a small number.

\end{proof}
\end{theorem}

We now prove the convergence rate of Algorithm \ref{alg:main}. The following lemma characterizes the relations between $\P(\bbb{x})$ and the optimal solution $\bbb{x}^*$ for any $\bbb{x}$. This is similar to the classical local proximal error bound in the literature \cite{luo1993error,TsengY09,Tseng10,yun2011block}.

\begin{lemma}\label{lemma:local:bound}
Assume $\bbb{x}$ is bounded. If $\bbb{x}$ is not the optimum of (\ref{eq:main}), there exists a constant $\eta\in(0,\infty)$ such that $\|\bbb{x} - \bbb{x}^* \| \leq \eta \|\bbb{x} - \P(\bbb{x})\|$.

\begin{proof}
\noi First, we prove that $\bbb{x} \neq \P(\bbb{x})$. This can be achieved by contradiction. According to the optimal condition of $\P(\bbb{x})$ in (\ref{eq:opt:bound0}), we obtain $\|\bbb{A}\P(\bbb{x}) + \bbb{b} +  \bbb{v}\| \leq  \|\bbb{C}\| \|\bbb{x}- \P(\bbb{x})\|$ with $\bbb{v}\in \partial h(\P(\bbb{x}))$. Assuming that $\bbb{x} = \P(\bbb{x})$, we obtain: $\bbb{A}\bbb{x} + \bbb{b} \in \partial h(\bbb{x})$, which contradicts with the condition that $\bbb{x}$ is not the optimal solution (refer to the optimality condition in (\ref{eq:main})). Therefore, it holds that $\bbb{x} \neq \P(\bbb{x})$. Second, by the boundedness of $\bbb{x}$ and $\bbb{x}^*$, there exists a sufficiently large constant $\eta\in(0,\infty)$ such that $\|\bbb{x} - \bbb{x}^* \| \leq \eta \|\bbb{x} - \P(\bbb{x})\|$.
\end{proof}
\end{lemma}


We now prove the convergence rate of Algorithm \ref{alg:main}.

\begin{theorem}
(Proof of Convergence Rate) We define $\delta \triangleq \tfrac{2\theta}{\omega} + \tfrac{2-\omega}{\omega}  \min(diag(\bbb{D}))$. Assume that $\omega$ and $\theta$ are chosen such that $\delta\in(0,\infty)$ and $\bbb{x}^k$ is bound for all $k$, we have:
\beq \label{eq:QQ}
\frac{f(\bbb{x}^{k+1}) - f(\bbb{x}^*)}{f(\bbb{x}^k)-f(\bbb{x}^*)} \leq \frac{C_1}{1+C_1}
\eeq
\noi where $C_1= ((3+\eta^2  \|\bbb{C}\|  + ( 2\eta^2 +2) \|\bbb{A}\|  )/{\delta}$. In other words, Algorithm \ref{alg:main} converges to the optimal solution Q-linearly.
\begin{proof}
Invoking Lemma \ref{lemma:local:bound} with $\bbb{x}=\bbb{x}^k$, we obtain:
\beq \label{eq:opt:bound:ineq}
\|\bbb{x}^k - \bbb{x}^*\| \leq \eta \|\bbb{x}^k - \P(\bbb{x}^k)\| ~\Rightarrow~\|\bbb{r}^k\| \leq \eta \|\bbb{d}^k\|
\eeq
\noi Invoking (\ref{eq:opt:bound0}) in Lemma \ref{lemma:opt:ineq} with $\bbb{x}=\bbb{x}^*,~\bbb{y}=\bbb{x}^k$, we derive the following inequalities:
\beq \label{eq:linear:conv}
& ~~~~~f^{k+1} - f^*\\
&\leq \textstyle  \la \bbb{r}^{k+1} ,\bbb{C} \bbb{d}^k  \ra   - \tfrac{1}{2}\la \bbb{r}^{k+1},\bbb{A}\bbb{r}^{k+1}\ra \\
&= \textstyle     \la \bbb{r}^{k}+\bbb{d}^{k} ,\bbb{C} \bbb{d}^k  \ra   - \tfrac{1}{2}\la \bbb{r}^{k}+\bbb{d}^{k},\bbb{A} (\bbb{r}^{k}+\bbb{d}^{k})\ra \\
&\leq \textstyle  \|\bbb{C}\|   (\|\bbb{d}^k\|\|\bbb{r}^{k}\| + \|\bbb{d}^k\|_2^2) + \tfrac{1}{2} \|\bbb{A}\| \|\bbb{r}^{k}+\bbb{d}^{k}\|_2^2 \\
&\leq \textstyle    \|\bbb{C}\|   (\tfrac{3}{2} \|\bbb{d}^k\|_2^2 + \tfrac{1}{2}\|\bbb{r}^{k}\|_2^2 )  + \|\bbb{A}\|(\|\bbb{r}^{k}\|_2^2 +\|\bbb{d}^{k}\|_2^2) \\
&\leq \textstyle      \|\bbb{C}\|\tfrac{3+\eta^2}{ 2} \|\bbb{d}^k\|_2^2  +  \|\bbb{A}\|( \eta^2 +1)\|\bbb{d}^k\|_2^2\\
&\leq \textstyle   ((3+\eta^2  \|\bbb{C}\|  + ( 2\eta^2 +2) \|\bbb{A}\|  )  \cdot (f^k-f^{k+1})/\delta \\
&= \textstyle C_1  (f^k-f^{k+1} ) \\
&= \textstyle C_1  (f^k-f^*) - C_1(f^{k+1}-f^*)
\eeq
\noi where the second step uses the fact that $\bbb{r}^{k+1}=\bbb{r}^{k}+\bbb{d}^{k}$; the third step uses the Cauchy-Schwarz inequality $\la \bbb{x},\bbb{y} \ra\leq \|\bbb{x}\|\|\bbb{y}\|,~\forall \bbb{x},\bbb{y}\in\mathbb{R}^n$ and the norm inequality $\|\bbb{Ax}\|\leq \|\bbb{A}\|\|\bbb{x}\|,~\forall \bbb{x} \in\mathbb{R}^n$; the fourth step uses the fact that ${1}/{2} \cdot \|\bbb{x}+\bbb{y}\|_2^2\leq \|\bbb{x}\|_2^2+\|\bbb{y}\|_2^2,~\forall \bbb{x},\bbb{y} \in\mathbb{R}^n$ and $ab\leq {1}/{2}\cdot a^2+{1}/{2} \cdot b^2,~\forall a,b \in \mathbb{R}$; the fifth step uses (\ref{eq:opt:bound:ineq}); the sixth step uses the descent condition in (\ref{eq:descent}). Rearranging the last inequality in (\ref{eq:linear:conv}), we have $(1+C_1)f(\bbb{x}^{k+1}) - f(\bbb{x}^*) \leq C_1 (f(\bbb{x}^k)-f(\bbb{x}^*))$ and obtain the inequality in (\ref{eq:QQ}). Since $\frac{C_1}{1+C_1}<1$, the sequence $\{f(\bbb{x}^k)\}_{k=0}^{\infty}$ converges to $f(\bbb{x}^*)$ linearly in the quotient sense.

\end{proof}
\end{theorem}

The following lemma is useful in our proof of iteration complexity.

\begin{lemma} \label{lemma:quadratic:recursive}
Suppose a nonnegative sequence $\{u^k\}_{k=0}^{\infty}$ satisfies $u^{k+1} \leq -2 C + 2C \sqrt{1+ \frac{u^k}{C} }$ for some constant $C\geq 0$. It holds that: $u^{k} \leq  \frac{C_2}{k}$, where $C_2=\max(8C,2\sqrt{Cu^0})$.
\begin{proof}
The proof of this lemma can be obtained by mathematical induction. (i) When $k=1$, we have $u^{1}  \leq -2 C + 2C \sqrt{1+ \frac{1}{C} u^0} \leq -2C + 2C (1+\sqrt{\frac{u^0}{C}} )=2\sqrt{Cu^0}\leq \frac{C_2}{k}$. (ii) When $k\geq 2$, we assume this result of this lemma is true. We derive the following results: $\frac{k+1}{k} \leq 2  \Rightarrow 4C \frac{k+1}{k} \leq 8C \leq C_2$$\Rightarrow  \frac{4C}{k(k+1)} \leq \frac{C_2}{(k+1)^2}$$\Rightarrow  4C \left(\frac{1}{k} - \frac{1}{k+1} \right)\leq \frac{C_2}{(k+1)^2}$$\Rightarrow   \frac{4C}{k} \leq    \frac{4C}{k+1} + \frac{C_2}{(k+1)^2}$$\Rightarrow    \frac{4CC_2}{k} \leq    \frac{4CC_2}{k+1} + \frac{C_2^2}{(k+1)^2}$$\Rightarrow   4C^2 \left( 1+ \frac{C_2}{k C}\right)  \leq  4C^2 +   \frac{4CC_2}{k+1} + \frac{C_2^2}{(k+1)^2}$$\Rightarrow   2C \sqrt{1+\frac{C_2}{Ck}} \leq 2C + \frac{C_2}{k+1}$$\Rightarrow   -2C + 2C \sqrt{1+\frac{C_2}{Ck}} \leq \frac{C_2}{k+1}$$\Rightarrow   u^{k+1} \leq \frac{C_2}{k+1}$.
\end{proof}
\end{lemma}

We now prove the iteration complexity of Algorithm \ref{alg:main}.

\begin{theorem}
(Proof of Iteration Complexity) We define $\delta \triangleq \tfrac{2\theta}{\omega} + \tfrac{2-\omega}{\omega}  \min(diag(\bbb{D}))$. Assume that $\omega$ and $\theta$ are chosen such that $\delta \in(0,\infty)$ and $\|\bbb{x}^k\|\leq R,~\forall k$, we have:
\beq
\textstyle
u^k  \leq
   \begin{cases}
   u^0(\frac{2C_4}{2C_4+1})^{k},  &\mbox{if~$\sqrt{f^{k}-f^{k+1}} \geq {C_3}/{C_4}$},~\forall k\leq \bar{k} \nn \\
   \frac{C_5}{k}, &\mbox{if~$\sqrt{f^{k}-f^{k+1}} < {C_3}/{C_4}$}, ~\forall k \geq  0
   \end{cases}
\eeq
\noi where $C_3=2R\|\bbb{C}\|(\tfrac{\delta}{2})^{1/2}$, $C_4=\tfrac{\delta}{2} (\|\bbb{C}\| + \|\bbb{A}\|  (\eta  + 1))$, $C_5=\max(8C_3^2,2C_3\sqrt{u^0})$, and $\bar{k}$ is some unknown iteration index.

\begin{proof}
Using similar strategies used in deriving (\ref{eq:linear:conv}), we have the following results:
\beq
&~~~~~~ u^{k+1}\\
&\leq \textstyle \la \bbb{r}^{k+1} ,\bbb{C} \bbb{d}^k  \ra   - \tfrac{1}{2}\la \bbb{r}^{k+1},\bbb{A}\bbb{r}^{k+1}\ra \\
&= \textstyle    \la \bbb{r}^{k}+\bbb{d}^{k} ,\bbb{C} \bbb{d}^k  \ra   - \tfrac{1}{2}\la \bbb{r}^{k}+\bbb{d}^{k},\bbb{A} (\bbb{r}^{k}+\bbb{d}^{k})\ra  \\
&\leq \textstyle    {\|\bbb{C}\|} (\|\bbb{r}^{k}\| \|\bbb{d}^k\| + \|\bbb{d}^{k}\|_2^2) + \tfrac{\|\bbb{A}\|}{2}\|\bbb{r}^{k}+\bbb{d}^{k}\|_2^2 \\
&\leq \textstyle    {\|\bbb{C}\|} (\|\bbb{r}^{k}\| \|\bbb{d}^k\| + \|\bbb{d}^{k}\|_2^2) + \|\bbb{A}\| (\|\bbb{r}^{k}\|_2^2+\|\bbb{d}^{k}\|_2^2 ) \\
&\leq \textstyle  {\|\bbb{C}\|} ( 2R\|\bbb{d}^k\| + \|\bbb{d}^{k}\|_2^2) + \|\bbb{A}\| (  \eta \| \bbb{d}^{k}\|_2^2+\|\bbb{d}^{k}\|_2^2 )\\
& \leq C_3 \sqrt{u^k-u^{k+1}} + C_4 (u^k-u^{k+1})\label{eq:recursion:u}
\eeq
\noi Now we consider the two cases for the recursion formula in (\ref{eq:recursion:u}): (i) $\sqrt{u^k-u^{k+1}} \geq \frac{C_3}{C_4}$ for some $k\leq \bar{k}$ (ii) $\sqrt{u^k-u^{k+1}} \leq \frac{C_3}{C_4}$ for all $k\geq 0$. In case (i), (\ref{eq:recursion:u}) implies that we have $u^{k+1}\leq  2 C_4 (u^{k}-u^{k+1})$ and rearranging terms gives: $u^{k+1}\leq \frac{2C_4}{2C_4+1} u^k$. Thus, we have: $u^{k+1}\leq (\frac{2C_4}{2C_4+1})^{k+1} u^0$. We now focus on case (ii). When $\sqrt{u^k-u^{k+1}} \leq \frac{C_3}{C_4}$, (\ref{eq:recursion:u}) implies that we have $u^{k+1}\leq  2 C_3 \sqrt{u^{k}-u^{k+1}}$ and rearranging terms yields:$\frac{(u^{k+1})^2}{4 C_3^2} + u^{k+1}  \leq   u^{k}$. Solving this quadratic inequality, we have: $u^{k+1} \leq  -2 C_3^2 + 2 C_3^2 \sqrt{1+\frac{1}{C_3^2} u^k}$; solving this recursive formulation by Lemma \ref{lemma:quadratic:recursive}, we obtain $ u^{k+1} \leq \frac{C_5}{k+1}$.

\end{proof}
\end{theorem}


Based on the discussions above, we have a few comments on Algorithm \ref{alg:main}. \textbf{(1)} When $h(\cdot)$ is empty and $\theta=0$, it reduces to the classical Gauss-Seidel method ($\omega=1$) and Successive Over-Relaxation method ($\omega\neq1$). \textbf{(2)} When $\bbb{A}$ contains zeros in its diagonal entries, one needs to set $\theta$ to a strictly positive number. This is to guarantee the strong convexity of the one dimensional subproblem and a bounded solution for any $h(\cdot)$. We remark that the introduction of the parameter $\theta$ is novel in this paper and it removes the assumption that $\bbb{A}$ is strictly positive-definite or strictly diagonally dominant, which is used in the classical result of GS and SOC method \cite{saad2003iterative,demmel1997applied}.

%% file: extension.tex
\begin{table*}[!t]
\scalebox{0.62}{\begin{tabular}{|c|c|c|c|c|c|c|c|}
  \hline
  data  & n & \cite{lin2007projected} & \cite{kim2011fast} &  \cite{kim2011fast}  & \cite{guan2012nenmf} & \cite{hsieh2011fast} & ours \\
    & & PG & AS &  BPP  & APG  & CGD & MSM \\
  \hline
  \multicolumn{8}{|c|}{\centering time limit=20} \\
20news &  20 &5.001e+06 & 2.762e+07 & 8.415e+06 & \cthree{4.528e+06} & \ctwo{4.515e+06} & \cone{4.506e+06} \\
20news &  50 &5.059e+06 & 2.762e+07 & 4.230e+07 & \cthree{3.775e+06} & \ctwo{3.544e+06} & \cone{3.467e+06} \\
20news &  100 &6.955e+06 & 5.779e+06 & 4.453e+07 & \ctwo{3.658e+06} & \cthree{3.971e+06} & \cone{2.902e+06} \\
20news &  200 &7.675e+06 & \ctwo{3.036e+06} & 1.023e+08 & \cthree{4.431e+06} & 3.573e+07 & \cone{2.819e+06} \\
20news &  300 &\cthree{1.997e+07} & 2.762e+07 & 1.956e+08 & \ctwo{4.519e+06} & 4.621e+07 & \cone{3.202e+06} \\
COIL &  20 &2.004e+09 & 5.480e+09 & 2.031e+09 & \cone{1.974e+09} & \cthree{1.976e+09} & \ctwo{1.975e+09} \\
COIL &  50 &1.412e+09 & 1.516e+10 & 6.962e+09 & \cthree{1.291e+09} & \ctwo{1.256e+09} & \cone{1.252e+09} \\
COIL &  100 &2.960e+09 & 2.834e+10 & 3.222e+10 & \cthree{9.919e+08} & \ctwo{8.745e+08} & \cone{8.510e+08} \\
COIL &  200 &3.371e+09 & 2.834e+10 & 5.229e+10 & \cthree{8.495e+08} & \ctwo{5.959e+08} & \cone{5.600e+08} \\
COIL &  300 &3.996e+09 & 2.834e+10 & 1.017e+11 & \cthree{8.493e+08} & \ctwo{5.002e+08} & \cone{4.956e+08} \\
TDT2 &  20 &1.597e+06 & 2.211e+06 & 1.688e+06 & \cone{1.591e+06} & \cthree{1.595e+06} & \ctwo{1.592e+06} \\
TDT2 &  50 &1.408e+06 & 2.211e+06 & 2.895e+06 & \cthree{1.393e+06} & \ctwo{1.390e+06} & \cone{1.385e+06} \\
TDT2 &  100 &1.300e+06 & 2.211e+06 & 6.187e+06 & \ctwo{1.222e+06} & \cthree{1.224e+06} & \cone{1.214e+06} \\
TDT2 &  200 &1.628e+06 & 2.211e+06 & 1.791e+07 & \ctwo{1.119e+06} & \cthree{1.227e+06} & \cone{1.079e+06} \\
TDT2 &  300 &1.915e+06 & \cthree{1.854e+06} & 3.412e+07 & \ctwo{1.172e+06} & 7.902e+06 & \cone{1.066e+06} \\
  \hline
\multicolumn{8}{|c|}{\centering time limit=30} \\
20news &  20 &4.716e+06 & 2.762e+07 & 7.471e+06 & \cthree{4.510e+06} & \ctwo{4.503e+06} & \cone{4.500e+06} \\
20news &  50 &4.569e+06 & 2.762e+07 & 5.034e+07 & \cthree{3.628e+06} & \ctwo{3.495e+06} & \cone{3.446e+06} \\
20news &  100 &6.639e+06 & 2.762e+07 & 4.316e+07 & \cthree{3.293e+06} & \ctwo{3.223e+06} & \cone{2.817e+06} \\
20news &  200 &\cthree{6.991e+06} & 2.762e+07 & 1.015e+08 & \ctwo{3.609e+06} & 7.676e+06 & \cone{2.507e+06} \\
20news &  300 &\cthree{1.354e+07} & 2.762e+07 & 1.942e+08 & \ctwo{4.519e+06} & 4.621e+07 & \cone{3.097e+06} \\
COIL &  20 &1.992e+09 & 4.405e+09 & 2.014e+09 & \cone{1.974e+09} & \cthree{1.975e+09} & \ctwo{1.975e+09} \\
COIL &  50 &1.335e+09 & 2.420e+10 & 5.772e+09 & \cthree{1.272e+09} & \ctwo{1.252e+09} & \cone{1.250e+09} \\
COIL &  100 &2.936e+09 & 2.834e+10 & 1.814e+10 & \cthree{9.422e+08} & \ctwo{8.623e+08} & \cone{8.458e+08} \\
COIL &  200 &3.362e+09 & 2.834e+10 & 4.627e+10 & \cthree{7.614e+08} & \ctwo{5.720e+08} & \cone{5.392e+08} \\
COIL &  300 &3.946e+09 & 2.834e+10 & 7.417e+10 & \cthree{6.734e+08} & \ctwo{4.609e+08} & \cone{4.544e+08} \\
TDT2 &  20 &1.595e+06 & 2.211e+06 & 1.667e+06 & \cone{1.591e+06} & \cthree{1.594e+06} & \ctwo{1.592e+06} \\
TDT2 &  50 &1.397e+06 & 2.211e+06 & 2.285e+06 & \cthree{1.393e+06} & \ctwo{1.389e+06} & \cone{1.385e+06} \\
TDT2 &  100 &1.241e+06 & 2.211e+06 & 5.702e+06 & \ctwo{1.216e+06} & \cthree{1.219e+06} & \cone{1.212e+06} \\
TDT2 &  200 &1.484e+06 & 1.878e+06 & 1.753e+07 & \ctwo{1.063e+06} & \cthree{1.104e+06} & \cone{1.049e+06} \\
TDT2 &  300 &1.879e+06 & 2.211e+06 & 3.398e+07 & \ctwo{1.060e+06} & \cthree{1.669e+06} & \cone{1.007e+06} \\
  \hline
\end{tabular}}\hspace{1pt}\scalebox{0.62}{\begin{tabular}{|c|c|c|c|c|c|c|c|}
  \hline
  data  &n & \cite{lin2007projected} & \cite{kim2011fast} &  \cite{kim2011fast}  & \cite{guan2012nenmf} & \cite{hsieh2011fast} & ours \\
    & & PG & AS &  BPP  & APG  & CGD & MSM \\
  \hline
\multicolumn{8}{|c|}{\centering time limit=40} \\
20news &  20 &4.622e+06 & 2.762e+07 & 7.547e+06 & \cone{4.495e+06} & \cthree{4.500e+06} & \ctwo{4.496e+06} \\
20news &  50 &4.386e+06 & 2.762e+07 & 1.562e+07 & \cthree{3.564e+06} & \ctwo{3.478e+06} & \cone{3.438e+06} \\
20news &  100 &6.486e+06 & 2.762e+07 & 4.223e+07 & \cthree{3.128e+06} & \ctwo{2.988e+06} & \cone{2.783e+06} \\
20news &  200 &6.731e+06 & 1.934e+07 & 1.003e+08 & \ctwo{3.304e+06} & \cthree{5.744e+06} & \cone{2.407e+06} \\
20news &  300 &\cthree{1.041e+07} & 2.762e+07 & 1.932e+08 & \ctwo{3.621e+06} & 4.621e+07 & \cone{2.543e+06} \\
COIL &  20 &1.987e+09 & 5.141e+09 & 2.010e+09 & \cone{1.974e+09} & \cthree{1.975e+09} & \ctwo{1.975e+09} \\
COIL &  50 &1.308e+09 & 2.403e+10 & 5.032e+09 & \cthree{1.262e+09} & \ctwo{1.250e+09} & \cone{1.248e+09} \\
COIL &  100 &2.922e+09 & 2.834e+10 & 2.086e+10 & \cthree{9.161e+08} & \ctwo{8.555e+08} & \cone{8.430e+08} \\
COIL &  200 &3.361e+09 & 2.834e+10 & 4.116e+10 & \cthree{7.075e+08} & \ctwo{5.584e+08} & \cone{5.289e+08} \\
COIL &  300 &3.920e+09 & 2.834e+10 & 7.040e+10 & \cthree{6.221e+08} & \ctwo{4.384e+08} & \cone{4.294e+08} \\
TDT2 &  20 &1.595e+06 & 2.211e+06 & 1.643e+06 & \cone{1.591e+06} & \cthree{1.594e+06} & \ctwo{1.592e+06} \\
TDT2 &  50 &1.394e+06 & 2.211e+06 & 1.933e+06 & \cthree{1.392e+06} & \ctwo{1.388e+06} & \cone{1.384e+06} \\
TDT2 &  100 &1.229e+06 & 2.211e+06 & 5.259e+06 & \ctwo{1.213e+06} & \cthree{1.216e+06} & \cone{1.211e+06} \\
TDT2 &  200 &1.389e+06 & 1.547e+06 & 1.716e+07 & \ctwo{1.046e+06} & \cthree{1.070e+06} & \cone{1.041e+06} \\
TDT2 &  300 &1.949e+06 & 1.836e+06 & 3.369e+07 & \ctwo{1.008e+06} & \cthree{1.155e+06} & \cone{9.776e+05} \\
  \hline
\multicolumn{8}{|c|}{\centering time limit=50} \\
20news &  20 &4.565e+06 & 2.762e+07 & 6.939e+06 & \cone{4.488e+06} & \cthree{4.498e+06} & \ctwo{4.494e+06} \\
20news &  50 &4.343e+06 & 2.762e+07 & 1.813e+07 & \cthree{3.525e+06} & \ctwo{3.469e+06} & \cone{3.432e+06} \\
20news &  100 &6.404e+06 & 2.762e+07 & 3.955e+07 & \cthree{3.046e+06} & \ctwo{2.878e+06} & \cone{2.765e+06} \\
20news &  200 &5.939e+06 & 2.762e+07 & 9.925e+07 & \ctwo{3.121e+06} & \cthree{4.538e+06} & \cone{2.359e+06} \\
20news &  300 &\cthree{9.258e+06} & 2.762e+07 & 1.912e+08 & \ctwo{3.621e+06} & 2.323e+07 & \cone{2.331e+06} \\
COIL &  20 &1.982e+09 & 7.136e+09 & 2.033e+09 & \cone{1.974e+09} & \cthree{1.975e+09} & \ctwo{1.975e+09} \\
COIL &  50 &1.298e+09 & 2.834e+10 & 4.365e+09 & \cthree{1.258e+09} & \ctwo{1.248e+09} & \cone{1.248e+09} \\
COIL &  100 &1.945e+09 & 2.834e+10 & 1.428e+10 & \cthree{9.014e+08} & \ctwo{8.516e+08} & \cone{8.414e+08} \\
COIL &  200 &3.362e+09 & 2.834e+10 & 3.760e+10 & \cthree{6.771e+08} & \ctwo{5.491e+08} & \cone{5.231e+08} \\
COIL &  300 &3.905e+09 & 2.834e+10 & 6.741e+10 & \cthree{5.805e+08} & \ctwo{4.226e+08} & \cone{4.127e+08} \\
TDT2 &  20 &1.595e+06 & 2.211e+06 & 1.622e+06 & \cone{1.591e+06} & \cthree{1.594e+06} & \ctwo{1.592e+06} \\
TDT2 &  50 &1.393e+06 & 2.211e+06 & 1.875e+06 & \cthree{1.392e+06} & \ctwo{1.386e+06} & \cone{1.384e+06} \\
TDT2 &  100 &1.223e+06 & 2.211e+06 & 4.831e+06 & \ctwo{1.212e+06} & \cthree{1.214e+06} & \cone{1.210e+06} \\
TDT2 &  200 &1.267e+06 & 2.211e+06 & 1.671e+07 & \ctwo{1.040e+06} & \cthree{1.054e+06} & \cone{1.036e+06} \\
TDT2 &  300 &1.903e+06 & 2.211e+06 & 3.328e+07 & \ctwo{9.775e+05} & \cthree{1.045e+06} & \cone{9.606e+05} \\
\hline
\end{tabular}}
\caption{Comparisons of objective values for non-negative matrix factorization for all the compared methods. The $1^{st}$, $2^{nd}$, and $3^{rd}$ best results are colored with \cone{red}, \ctwo{blue} and \cthree{green}, respectively.}
\label{tab:nmf}
\end{table*}

\section{Extensions}\label{sect:extension}
This section discusses several extensions of our proposed matrix splitting method for solving (\ref{eq:main}).

\subsection{When h is Nonconvex}
When $h(\bbb{x})$ is nonconvex, our theoretical analysis breaks down in (\ref{eq:fail}) and the exact solution to the triangle proximal operator $\P(\bbb{x}^k)$ in (\ref{eq:subproblem}) cannot be guaranteed. However, our Gaussian elimination procedure in Algorithm \ref{alg:sub} can still be applied. What one needs is to solve a one-dimensional nonconvex subproblem in (\ref{eq:1d:subp}). For example, when $h(t)=\lambda |t|_0$ (\eg in the case of the $\ell_0$ norm), it has an analytical solution: $t^* = {\tiny \left\{
        \begin{array}{cc}
          -\bbb{w}_j/\bbb{B}_{j,j}, & {\bbb{w}_j^2 > 2\lambda  \bbb{B}_{j,j}} \\
          0, & {\bbb{w}_j^2 \leq 2\lambda  \bbb{B}_{j,j}}
        \end{array}
      \right.}$; when $h(t)=\lambda |{t}|^p$ and $p<1$, it admits a closed form solution for some special values \cite{xu2012regularization}, such as $p=\frac{1}{2}$ or $\frac{2}{3}$.

Our matrix splitting method is guaranteed to converge even when ${h}(\cdot)$ is nonconvex. Specifically, we present the following theorem.

\begin{theorem} 
(Proof of Global Convergence when ${h}(\cdot)$ is Nonconvex) Assume the nonconvex one-dimensional subproblem in (\ref{eq:1d:subp}) can be solved globally and analytically. We define $\delta\triangleq \min\left(\theta/\omega +(1-\omega)/\omega \cdot diag(\bbb{D})\right)$. If we choose $\omega$ and $\theta$ such that $\delta\in(0,\infty)$, we have: (i)
\beq \label{eq:nonconvex:suf:dec}
\textstyle f(\bbb{x}^{k+1}) - f(\bbb{x}^k) \leq - \frac{\delta}{2} \|\bbb{x}^{k+1}-\bbb{x}^k\|_2^2 \leq 0
\eeq
\noi (ii) Algorithm \ref{alg:main} is globally convergent.
\begin{proof}
(i) Due to the optimality of the one-dimensional subproblem in (\ref{eq:1d:subp}), for all $j=1,2,...,n$, we have:
\beq
&\textstyle  \tfrac{1}{2}\bbb{B}_{j,j}(\bbb{x}^{k+1}_j)^2 + (\bbb{u}_j +  \sum_{i=1}^{j-1} \bbb{B}_{j,i}\bbb{x}^{k+1}_{i}) \bbb{x}^{k+1}_j +  h(\bbb{x}^{k+1}_j) \nn \\
&\leq  \textstyle\tfrac{1}{2}\bbb{B}_{j,j}{t}_j^2 + (\bbb{u}_j +  \sum_{i=1}^{j-1} \bbb{B}_{j,i}\bbb{x}^{k+1}_{i}) {t}_j +  h({t}_j),~\forall {t}_j
\eeq
\noi Letting ${t}_1=\bbb{x}^k_1,~{t}_2=\bbb{x}^k_2,~...~,{t}_n=\bbb{x}^k_n$, we obtain:
\beq
&\textstyle \tfrac{1}{2} \sum_{i}^n\bbb{B}_{i,i} (\bbb{x}^{k+1})^2 + \la \bbb{u} + \bbb{Lx}^{k+1},\bbb{x}^{k+1} \ra  +  h(\bbb{x}^{k+1})\nn\\
&\textstyle \leq\tfrac{1}{2}\sum_{i}^n\bbb{B}_{i,i} (\bbb{x}^{k})^2 + \la \bbb{u}+\bbb{Lx}^{k+1},\bbb{x}^{k} \ra +  h(\bbb{x}^{k})\nn\\
\eeq
\noi Since $\bbb{u} = \bbb{b}+\bbb{C}\bbb{x}^k$, we obtain the following inequality:
\beq
\textstyle f^{k+1} + \tfrac{1}{2}\la \bbb{x}^{k+1},(\tfrac{1}{\omega}(\bbb{D}+\theta \bbb{I})+2\bbb{L}-\bbb{A})\bbb{x}^{k+1} + 2 \bbb{C}\bbb{x}^k \ra\nn\\
\leq   f^k + \tfrac{1}{2} \la \bbb{x}^{k},(\tfrac{1}{\omega}(\bbb{D}+\theta \bbb{I})+2\bbb{C}-\bbb{A}) \bbb{x}^{k} + 2\bbb{L}\bbb{x}^{k+1}\ra~~~~\nn
\eeq
\noi By denoting $\bbb{S} \triangleq \bbb{L}-\bbb{L}^T$ and $\bbb{T}\triangleq((\omega-1)\bbb{D}-\theta \bbb{I})/\omega$, we have: $\tfrac{1}{\omega}(\bbb{D}+\theta \bbb{I})+2\bbb{L}-\bbb{A}=\bbb{T}-\bbb{S}$, $\tfrac{1}{\omega}(\bbb{D}+\theta \bbb{I})+2\bbb{C}-\bbb{A}=\bbb{S}-\bbb{T}$, and $\bbb{L}-\bbb{C}^T=-\bbb{T}$. Therefore, we have the following inequalities:
\beq
&\textstyle  f^{k+1} -  f^k\leq  \tfrac{1}{2}\la \bbb{x}^{k+1},(\bbb{T} - \bbb{S})\bbb{x}^{k+1}\ra   - \la \bbb{x}^k, \bbb{T}\bbb{x}^{k+1})   \nn\\
& \textstyle + \tfrac{1}{2}\la \bbb{x}^k,(\bbb{T} - \bbb{S})\bbb{x}^k\ra  =   \tfrac{1}{2}\la \bbb{x}^k-\bbb{x}^{k+1}, \bbb{T}(\bbb{x}^k-\bbb{x}^{k+1})\ra \nn\\
&\leq -\tfrac{\delta }{2}\|\bbb{x}^{k+1}-\bbb{x}^{k}\|_2^2
\eeq
\noi where the first equality uses $\la \bbb{x},\bbb{S}\bbb{x}\ra =0~\forall \bbb{x}$, since $\bbb{S}$ is a Skew-Hermitian matrix. The last step uses $\bbb{T} + \delta \bbb{I} \preceq \bbb{0}$, since $\bbb{x}+ \min(-\bbb{x})\leq \bbb{0}~\forall \bbb{x}$. Thus, we obtain the sufficient decrease inequality in (\ref{eq:nonconvex:suf:dec}). 

(ii) Based on the sufficient decrease inequality in (\ref{eq:nonconvex:suf:dec}), we have: $f(\bbb{x}^k)$ is a non-increasing sequence, $\|\bbb{x}^k-\bbb{x}^{k+1}\|\rightarrow 0$, and $f(\bbb{x}^{k+1})<f(\bbb{x}^k)$ if $\bbb{x}^k\neq\bbb{x}^{k+1}$. We note that (\ref{eq:opt:bound0}) can be still applied even ${h}(\cdot)$ is nonconvex. Using the same methodology as in the second part of Theorem \ref{theorem:1}, we obtain that $\partial f(\bbb{x}^k)\rightarrow \bbb{0}$, which implies the convergence of the algorithm.

Note that guaranteeing $\delta\in(0,\infty)$ can be achieved by simply choosing $\omega\in(0,1)$ and setting $\theta$ to a small number.

\end{proof}
\end{theorem}


\subsection{When x is a Matrix}
In many applications (e.g. nonegative matrix factorization and sparse coding), the solutions exist in the matrix form as follows: $\min_{\bbb{X}\in \mathbb{R}^{n\times r}}~~\tfrac{1}{2}tr(\bbb{X}^T\bbb{A}\bbb{X}) + tr(\bbb{X}^T\bbb{R})+ h(\bbb{X})$, where $\bbb{R}\in \mathbb{R}^{n\times r}$. Our matrix splitting algorithm can still be applied in this case. Using the same technique to decompose $\bbb{A}$ as in (\ref{eq:matrix}): $\bbb{A}=\bbb{B}+\bbb{C}$, one needs to replace (\ref{eq:subproblem}) to solve the following nonlinear equation: $\bbb{BZ}^* + \bbb{U} +  \partial h(\bbb{Z}^*) \in 0$, where $\bbb{U}=\bbb{R}+\bbb{C}\bbb{X}^k$. It can be decomposed into $r$ independent components. By updating every column of $\bbb{X}$, the proposed algorithm can be used to solve the matrix problem above. Thus, our algorithm can also make good use of existing parallel architectures to solve the matrix optimization problem.



\subsection{When q is not Quadratic}
Following previous work \cite{TsengY09,yuan2014newton}, one can approximate the objective around the current solution $\bbb{x}^k$ by the second-order Taylor expansion of $q(\cdot)$: $Q(\bbb{y},\bbb{x}^k)\triangleq q(\bbb{\bbb{x}^k}) + \la \nabla q(\bbb{\bbb{x}^k}) , \bbb{y} - \bbb{\bbb{x}^k}\ra + \frac{1}{2} (\bbb{y} - \bbb{\bbb{x}^k})^T \nabla^2 q(\bbb{\bbb{x}^k}) (\bbb{y} - \bbb{\bbb{x}^k}) + h(\bbb{y})$, where $\nabla q(\bbb{\bbb{x}^k})$ and $\nabla^2 q(\bbb{\bbb{x}^k})$ denote the gradient and Hessian of $q(\bbb{x})$ at $\bbb{x}^k$, respectively. In order to generate the next solution that decreases the objective, one can minimize the quadratic model above by solving: $\bar{\bbb{x}}^k = \arg\min_{\bbb{y}}~Q(\bbb{y},\bbb{x}^k)$. The new estimate is obtained by the update: $\bbb{x}^{k+1}=\bbb{x}^k + \beta (\bar{\bbb{x}}^k-\bbb{x}^k)$, where $\beta\in(0,1)$ is the step-size selected by backtracking line search.

%% file: experiment.tex
\begin{figure*} 
\centering

      \begin{subfigure}{\fivefigwid}\includegraphics[width=\textwidth,height=\objimghei]{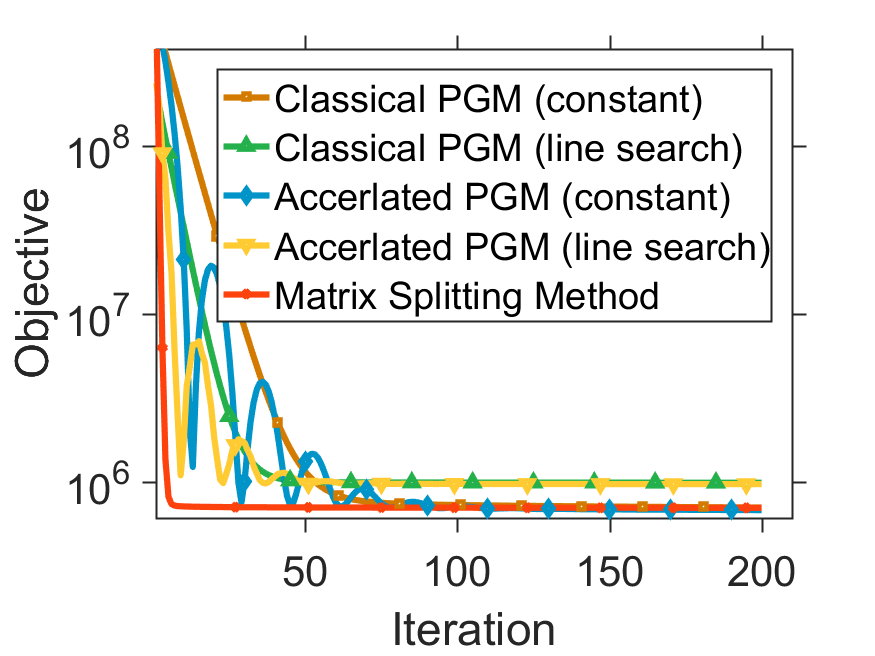}\vspace{-6pt} \caption{\footnotesize $\lambda=5$}\end{subfigure}\ghs
      \begin{subfigure}{\fivefigwid}\includegraphics[width=\textwidth,height=\objimghei]{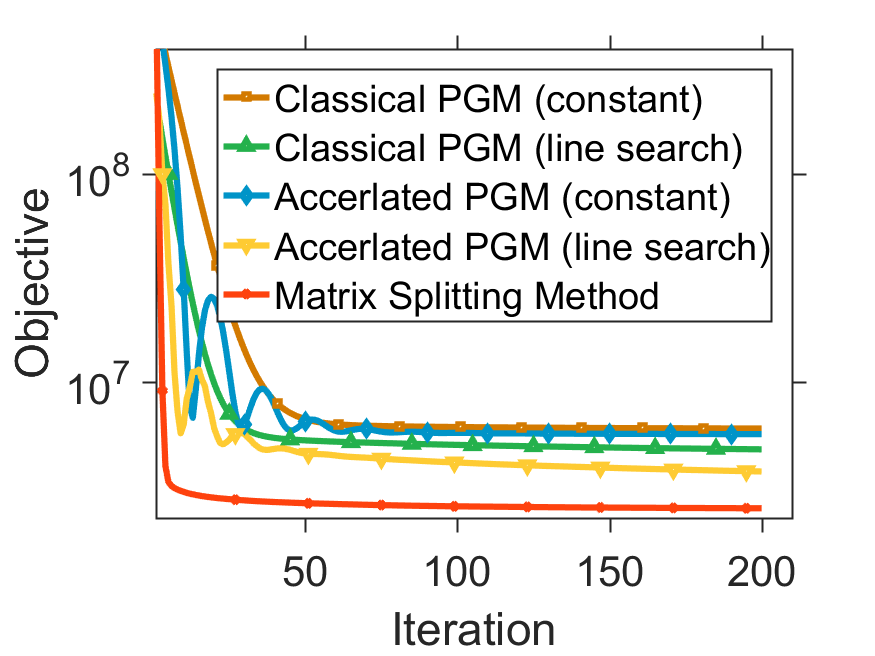}\vspace{-6pt} \caption{\footnotesize $\lambda=50$}\end{subfigure}\ghs
      \begin{subfigure}{\fivefigwid}\includegraphics[width=\textwidth,height=\objimghei]{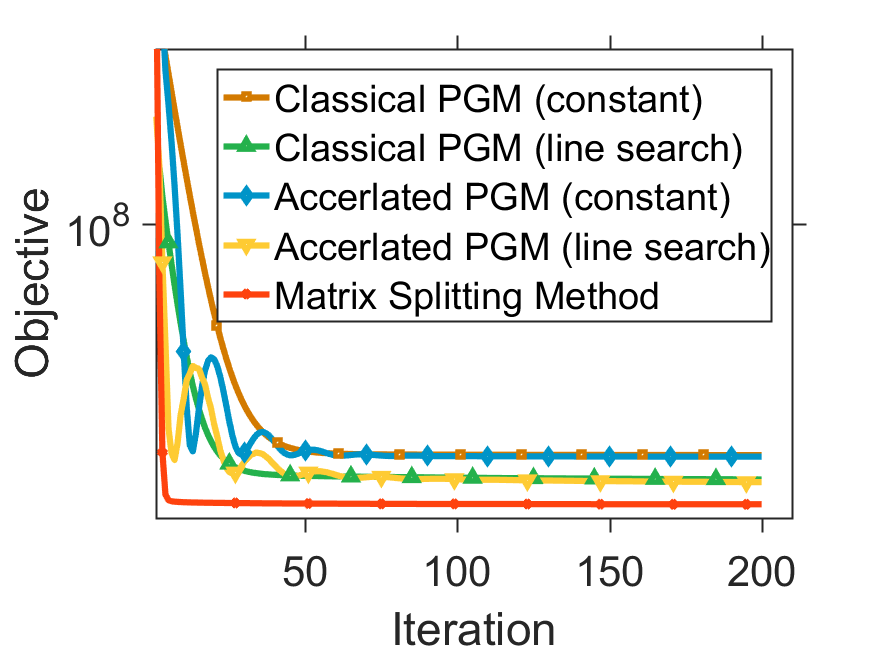}\vspace{-6pt} \caption{\footnotesize $\lambda=500$}\end{subfigure}\ghs
      \begin{subfigure}{\fivefigwid}\includegraphics[width=\textwidth,height=\objimghei]{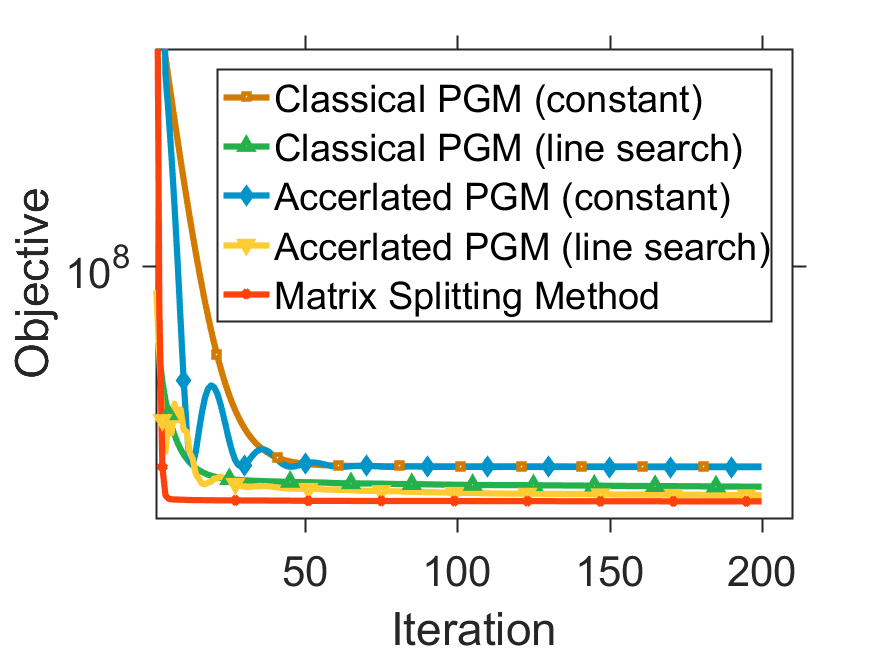}\vspace{-6pt} \caption{\footnotesize $\lambda=5000$}\end{subfigure}\ghs
      \begin{subfigure}{\fivefigwid}\includegraphics[width=\textwidth,height=\objimghei]{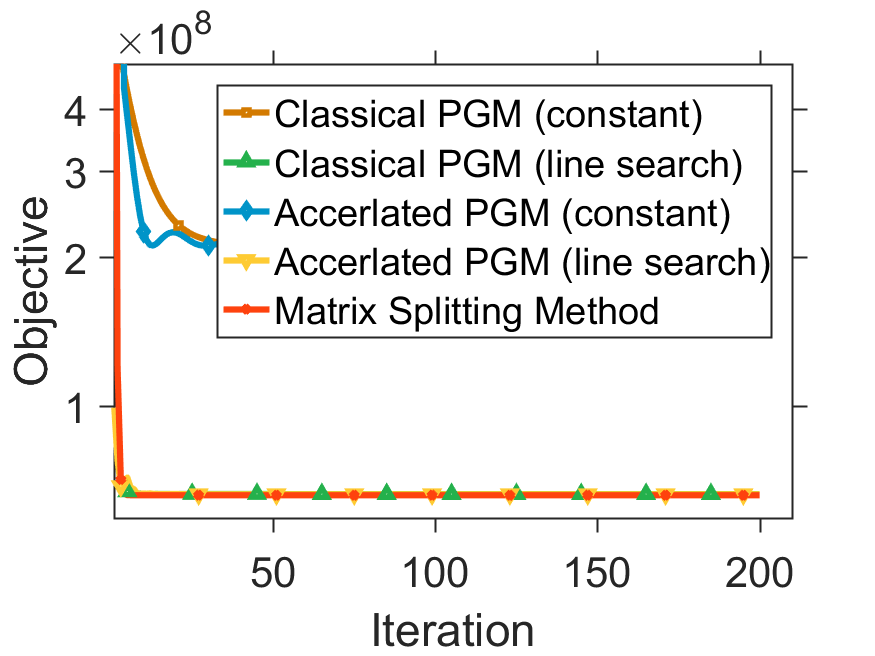}\vspace{-6pt} \caption{\footnotesize $\lambda=50000$}\end{subfigure}\\

      \begin{subfigure}{\fivefigwid}\includegraphics[width=\textwidth,height=\objimghei]{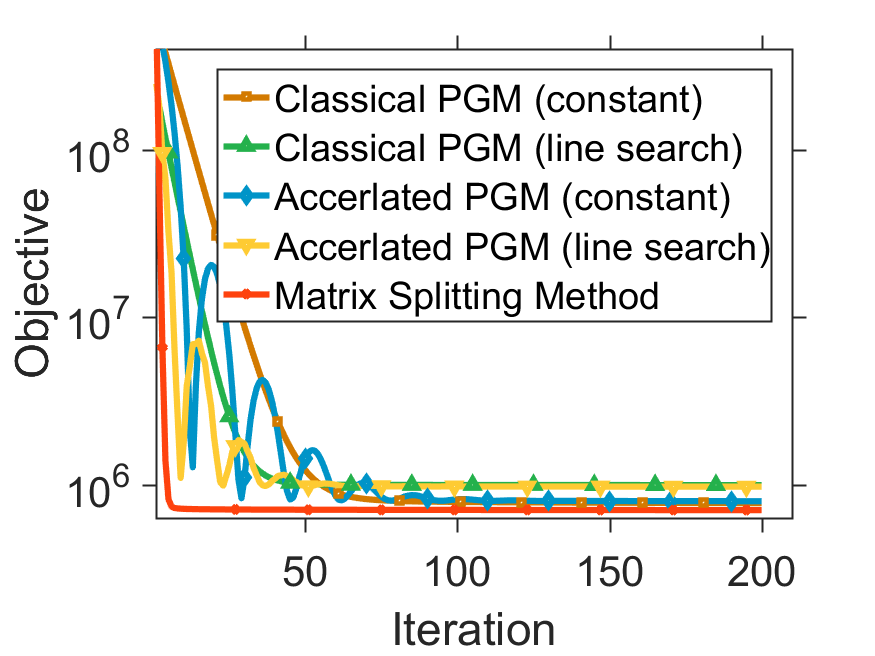}\vspace{-6pt} \caption{\footnotesize $\lambda=5$}\end{subfigure}\ghs
      \begin{subfigure}{\fivefigwid}\includegraphics[width=\textwidth,height=\objimghei]{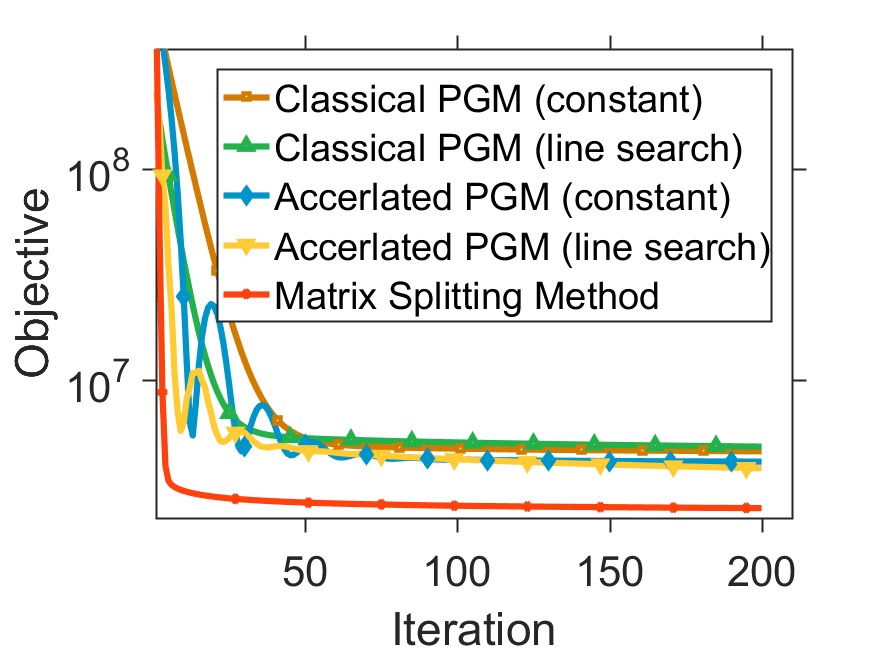}\vspace{-6pt} \caption{\footnotesize $\lambda=50$}\end{subfigure}\ghs
      \begin{subfigure}{\fivefigwid}\includegraphics[width=\textwidth,height=\objimghei]{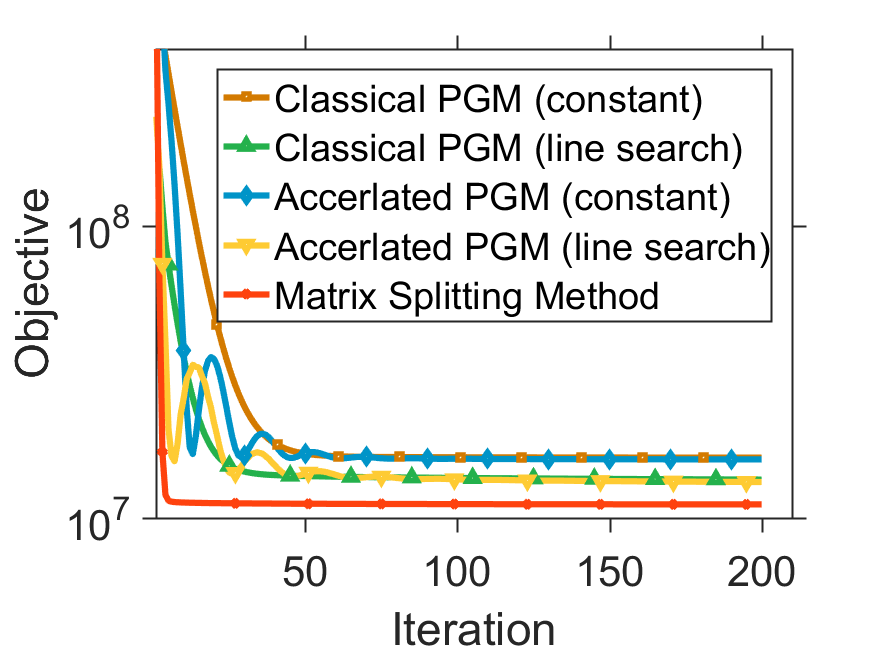}\vspace{-6pt} \caption{\footnotesize $\lambda=500$}\end{subfigure}\ghs
      \begin{subfigure}{\fivefigwid}\includegraphics[width=\textwidth,height=\objimghei]{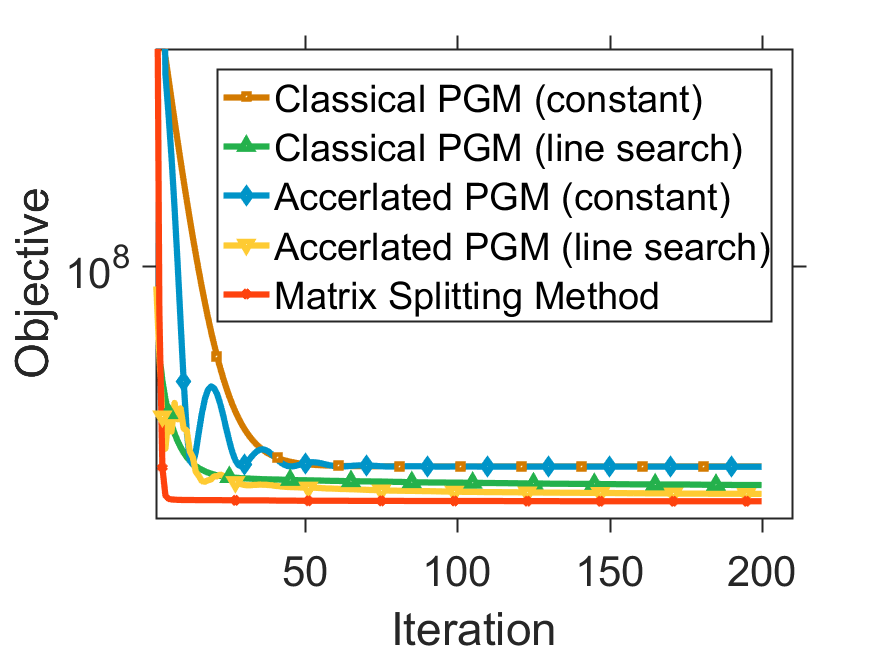}\vspace{-6pt} \caption{\footnotesize $\lambda=5000$}\end{subfigure}\ghs
      \begin{subfigure}{\fivefigwid}\includegraphics[width=\textwidth,height=\objimghei]{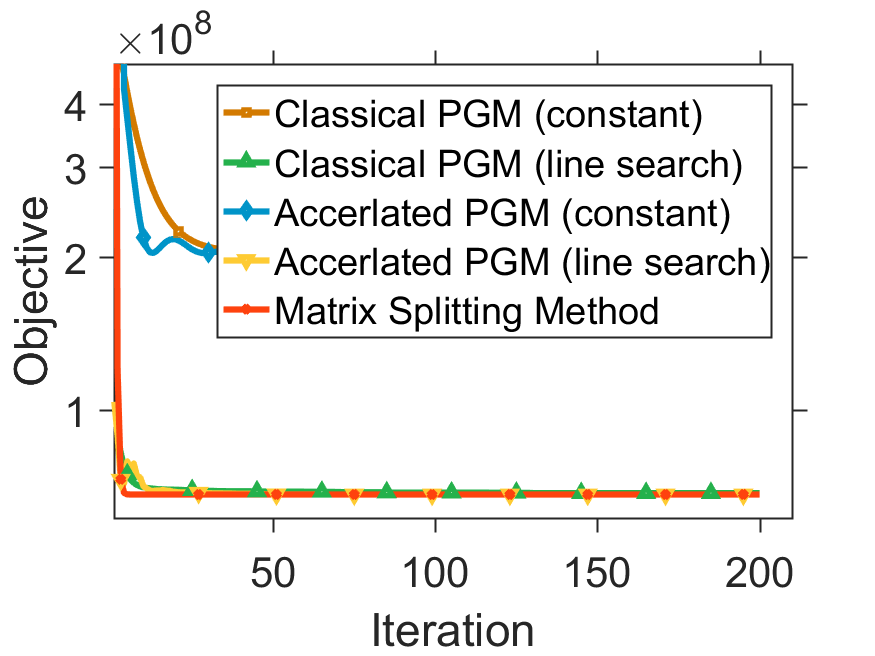}\vspace{-6pt} \caption{\footnotesize $\lambda=50000$}\end{subfigure}\\

      \begin{subfigure}{\fivefigwid}\includegraphics[width=\textwidth,height=\objimghei]{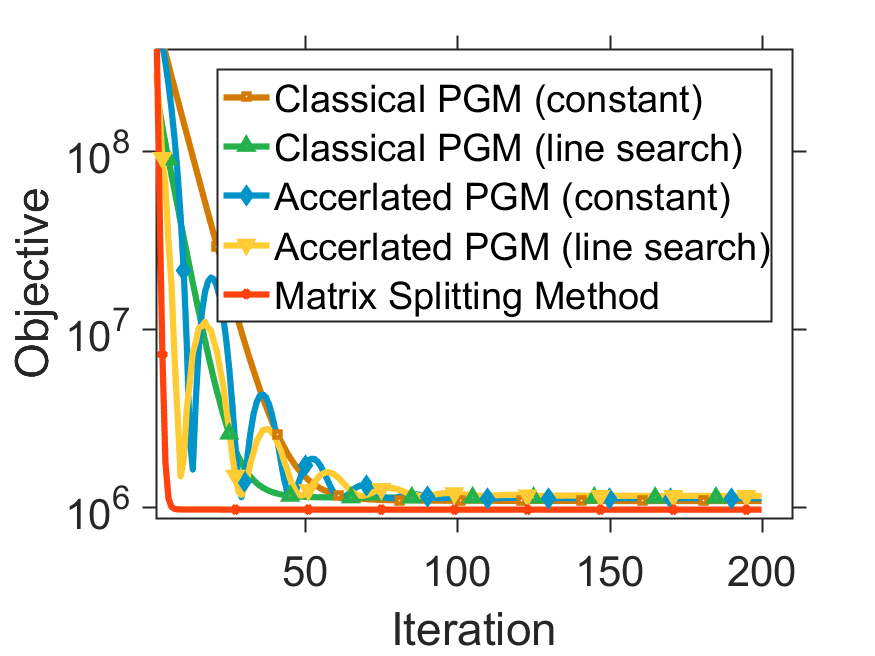}\vspace{-6pt} \caption{\footnotesize $\lambda=5$}\end{subfigure}\ghs
      \begin{subfigure}{\fivefigwid}\includegraphics[width=\textwidth,height=\objimghei]{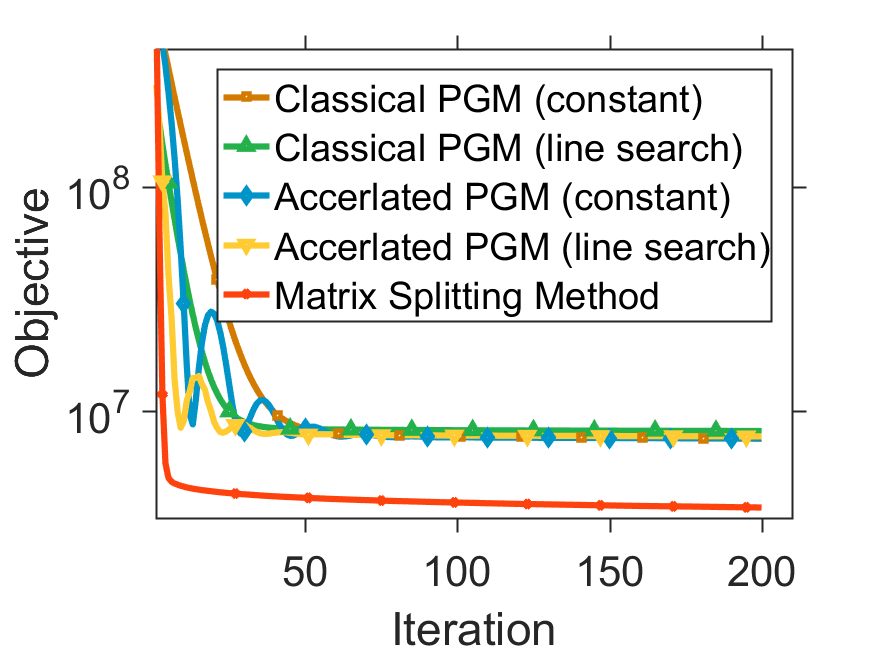}\vspace{-6pt} \caption{\footnotesize $\lambda=50$}\end{subfigure}\ghs
      \begin{subfigure}{\fivefigwid}\includegraphics[width=\textwidth,height=\objimghei]{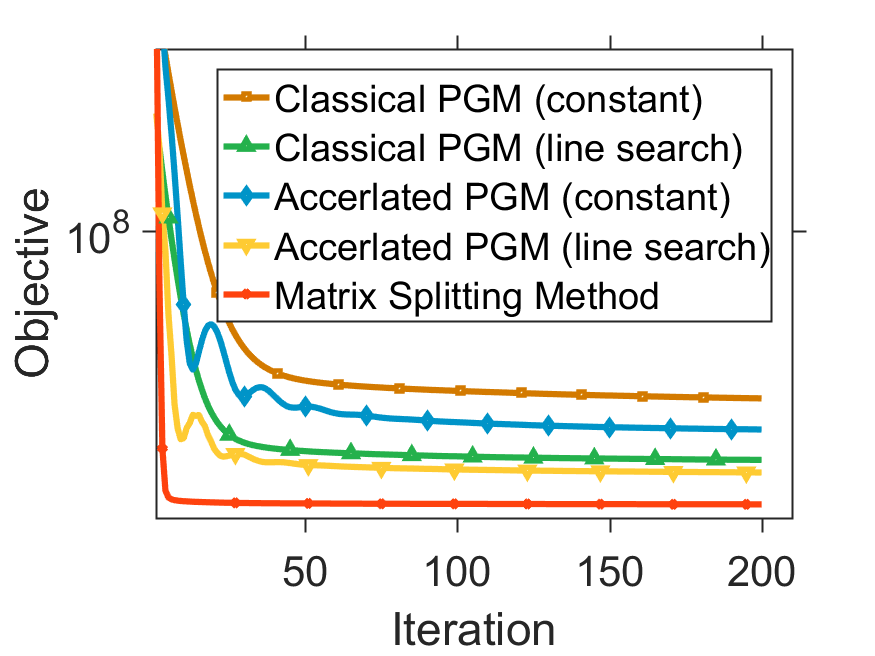}\vspace{-6pt} \caption{\footnotesize $\lambda=500$}\end{subfigure}\ghs
      \begin{subfigure}{\fivefigwid}\includegraphics[width=\textwidth,height=\objimghei]{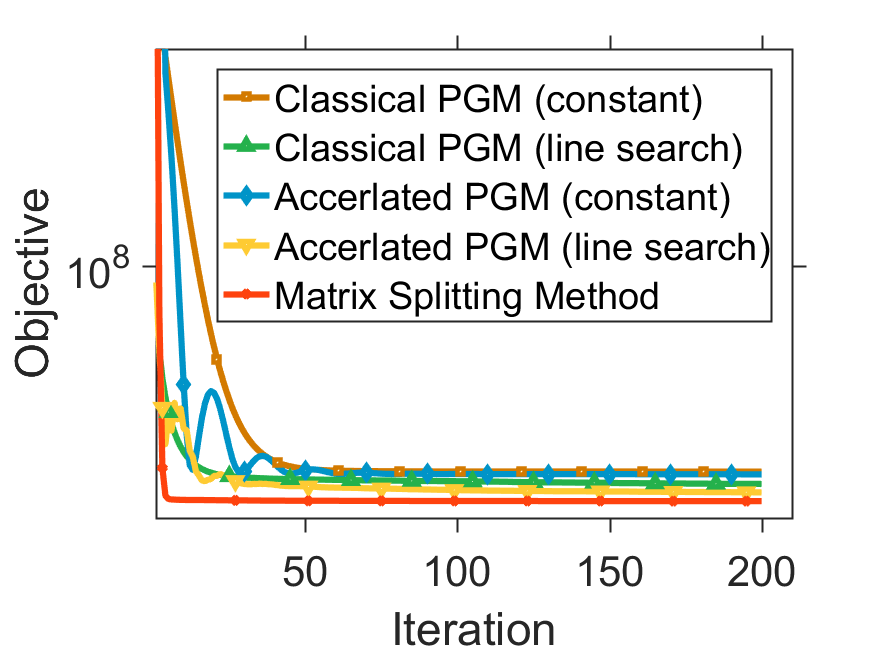}\vspace{-6pt} \caption{\footnotesize $\lambda=5000$}\end{subfigure}\ghs
      \begin{subfigure}{\fivefigwid}\includegraphics[width=\textwidth,height=\objimghei]{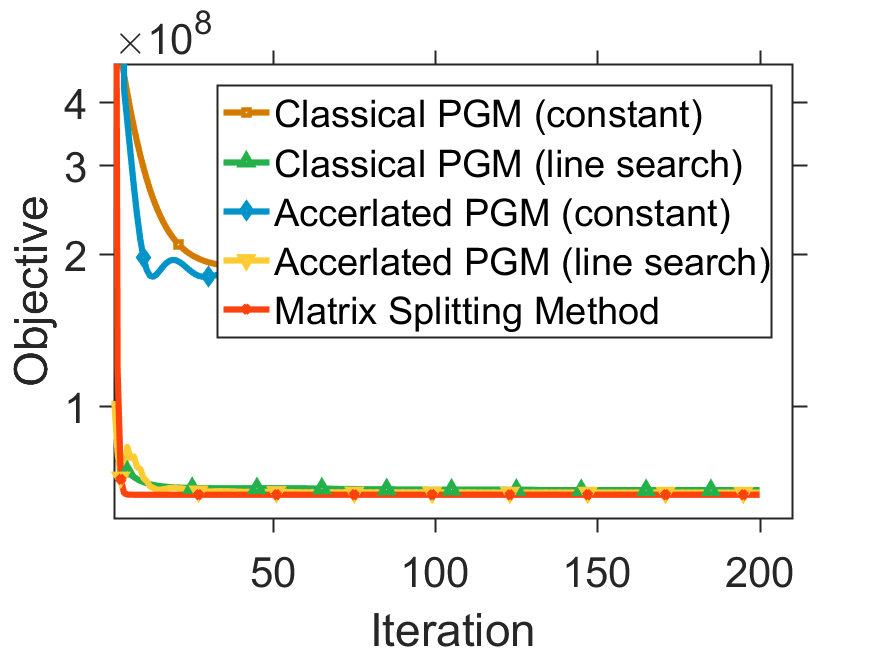}\vspace{-6pt} \caption{\footnotesize $\lambda=50000$}\end{subfigure}\\

      \begin{subfigure}{\fivefigwid}\includegraphics[width=\textwidth,height=\objimghei]{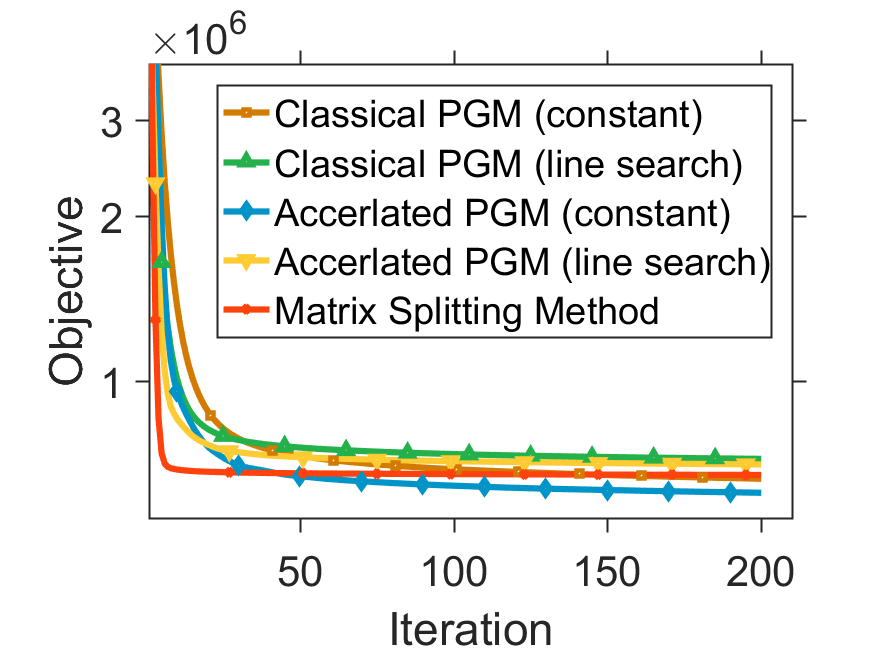}\vspace{-6pt} \caption{\footnotesize $\lambda=5$}\end{subfigure}\ghs
      \begin{subfigure}{\fivefigwid}\includegraphics[width=\textwidth,height=\objimghei]{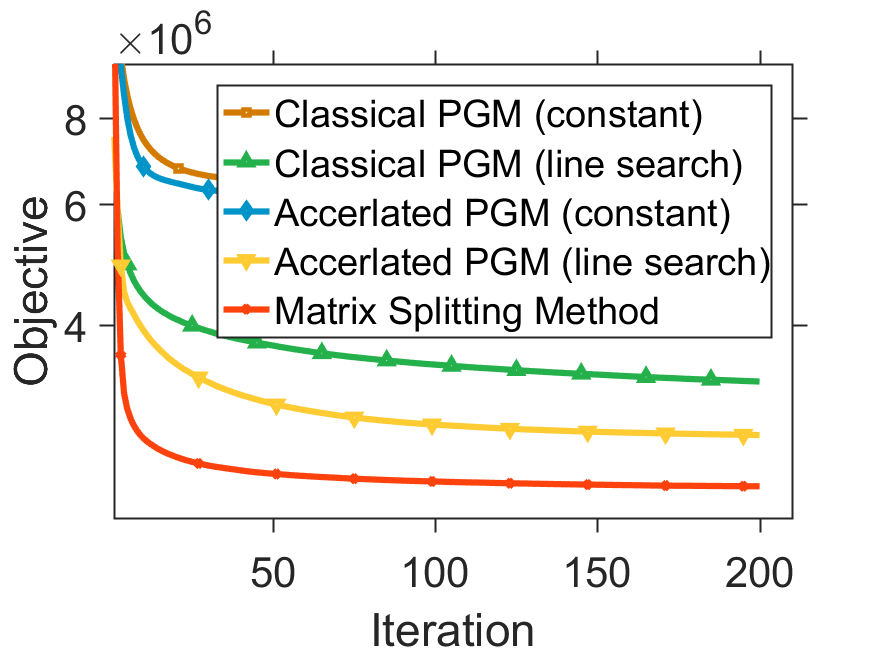}\vspace{-6pt} \caption{\footnotesize $\lambda=50$}\end{subfigure}\ghs
      \begin{subfigure}{\fivefigwid}\includegraphics[width=\textwidth,height=\objimghei]{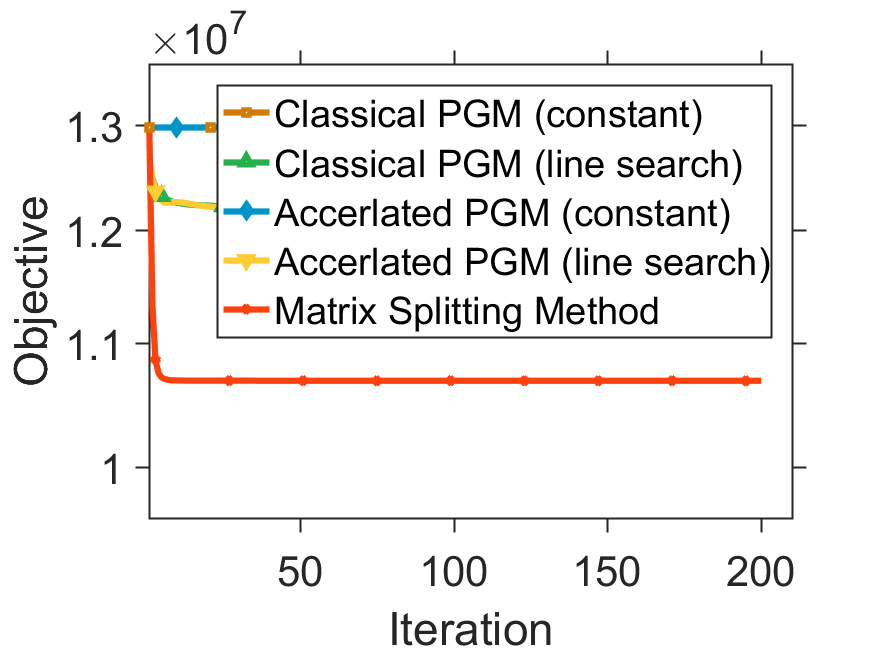}\vspace{-6pt} \caption{\footnotesize $\lambda=500$}\end{subfigure}\ghs
      \begin{subfigure}{\fivefigwid}\includegraphics[width=\textwidth,height=\objimghei]{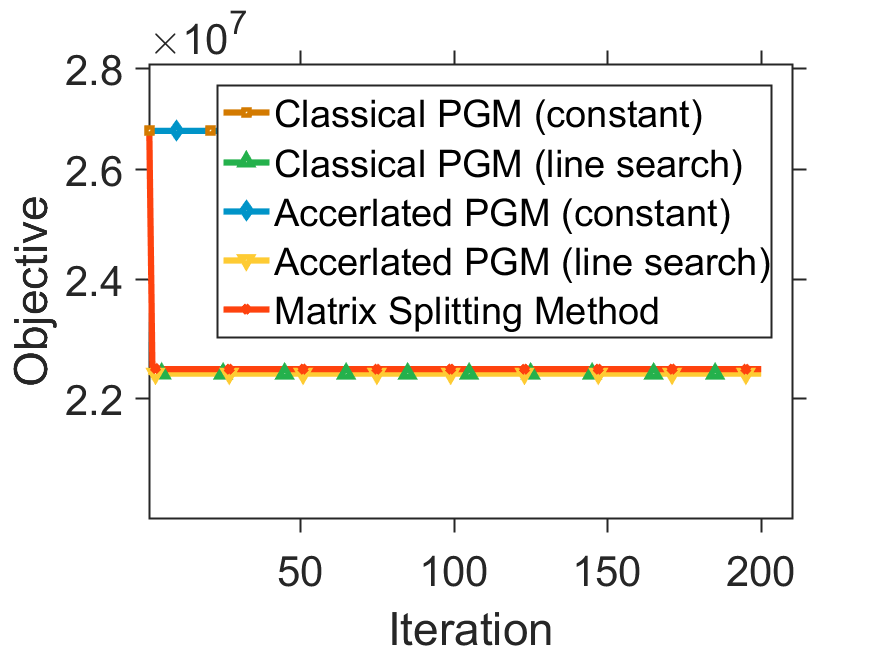}\vspace{-6pt} \caption{\footnotesize $\lambda=5000$}\end{subfigure}\ghs
      \begin{subfigure}{\fivefigwid}\includegraphics[width=\textwidth,height=\objimghei]{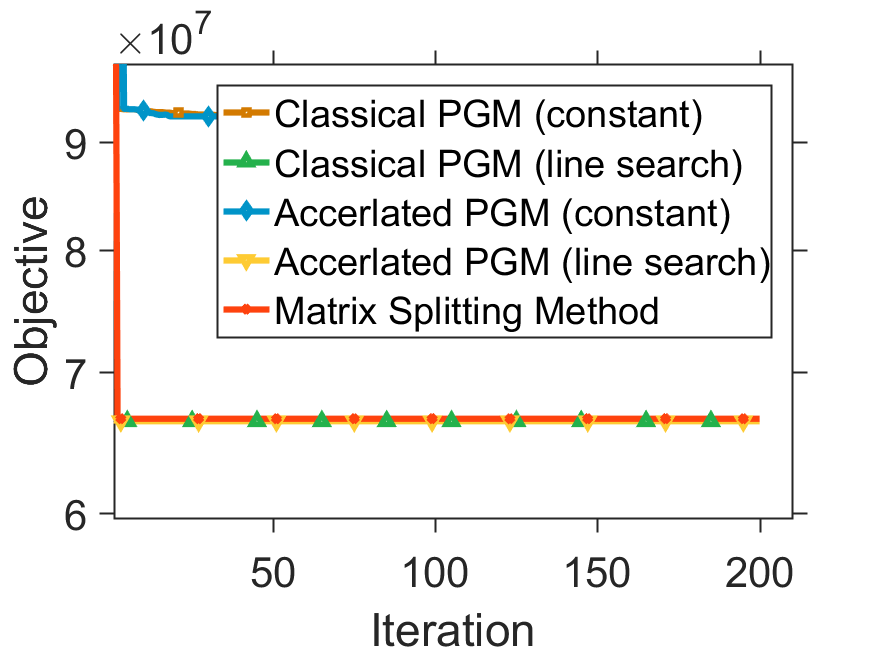}\vspace{-6pt} \caption{\footnotesize $\lambda=50000$}\end{subfigure}\\
\caption{Convergence behavior for solving (\ref{eq:card:sparse:coding}) with fixing $\bbb{W}$ for different $\lambda$ and initializations. Denoting $\tilde{\bbb{O}}$ as an arbitrary standard Gaussian random matrix of suitable size, we consider the following four initializations for $\bbb{H}$. First row: $\bbb{H}=0.1 \times \tilde{\bbb{O}}$. Second row: $\bbb{H}=1\times \tilde{\bbb{O}}$. Third row: $\bbb{H}=10\times \tilde{\bbb{O}}$. Fourth row: $\bbb{H}$ is set to the output of the orthogonal matching pursuit.}
\label{fig:convergece:obj}
\end{figure*}

\section{Experiments}

This section demonstrates the efficiency and efficacy of the proposed Matrix Splitting Method (MSM) by considering two important applications: nonnegative matrix factorization (NMF) \cite{lee1999learning,lin2007projected} and cardinality regularized sparse coding \cite{olshausen1996emergence,Quan2016CVPR,lee2006efficient}. We implement our method in MATLAB on an Intel 2.6 GHz CPU with 8 GB RAM. Only our generalized Gaussian elimination procedure is developed in C and wrapped into the MATLAB code, since it requires an elementwise loop that is quite inefficient in native MATLAB. We report our results with the choice $\theta=0.01$ and $\omega=1$ in all our experiments. 

\subsection{Nonnegative Matrix Factorization (NMF)}

Nonnegative matrix factorization \cite{lee1999learning} is a very useful tool for feature extraction and identification in the fields of text mining and image understanding. It is formulated as the following optimization problem:
\beq
\textstyle \underset{\bbb{W},\bbb{H}}\min~~\frac{1}{2}\|\bbb{Y}-\bbb{WH}\|_F^2 ,~~s.t.~~\bbb{W}\geq 0,~\bbb{H}\geq 0
\eeq
\noi where $\bbb{W}\in\mathbb{R}^{m\times n}$ and $\bbb{H}\in\mathbb{R}^{n\times d}$. Following previous work \cite{kim2011fast,guan2012nenmf,lin2007projected,hsieh2011fast}, we alternatively minimize the objective while keeping one of the two variables fixed. In each alternating subproblem, we solve a convex nonnegative least squares problem, where our MSM method is used. We conduct experiments on three datasets \footnote{\url{http://www.cad.zju.edu.cn/home/dengcai/Data/TextData.html}} 20news, COIL, and TDT2. The size of the datasets are $18774\times 61188,~7200\times 1024,~9394\times 36771$, respectively. We compare MSM against the following state-of-the-art methods: (1) Projective Gradient (PG)  \cite{lin2007projected,bertsekas1999nonlinear} that updates the current solution via steep gradient descent and then maps a point back to the bounded feasible region \footnote{\url{https://www.csie.ntu.edu.tw/~cjlin/libmf/}}; (2) Active Set (AS) method \cite{kim2011fast}; (3) Block Principal Pivoting (BPP) method \cite{kim2011fast} \footnote{\url{http://www.cc.gatech.edu/~hpark/nmfsoftware.php}} that iteratively identifies an active and passive set by a principal pivoting procedure and solves a reduced linear system; (4) Accelerated Proximal Gradient (APG) \cite{guan2012nenmf} \footnote{\url{https://sites.google.com/site/nmfsolvers/}} that applies Nesterov's momentum strategy with a constant step size to solve the convex sub-problems; (5) Coordinate Gradient Descent (CGD) \cite{hsieh2011fast} \footnote{\url{http://www.cs.utexas.edu/~cjhsieh/nmf/}} that greedily selects one coordinate by measuring the objective reduction and optimizes for a single variable via closed-form update. Similar to our method, the core procedure of CGD is developed in C and wrapped into the MATLAB code, while all other methods are implemented using builtin MATLAB functions.

We use the same settings as in \cite{lin2007projected}. We compare objective values after running $t$ seconds with $t$ varying from 20 to 50. Table \ref{tab:nmf} presents average results of using 10 random initial points, which are generated from a standard normal distribution. While the other methods may quickly lower objective values when $n$ is small ($n=20$), MSM catches up very quickly and achieves a faster convergence speed when $n$ is large. It generally achieves the best performance in terms of objective value among all the methods.

\def\boxscale{0.66}
\begin{table*}
\centering
\scalebox{\boxscale}{
\begin{tabular}{|p{2.35cm}|p{0.8cm}|p{0.8cm}|p{0.8cm}|p{0.8cm}|p{0.8cm}|}
\hline
 \multicolumn{2}{|c|}{} & \multicolumn{2}{c|}{OMP init} & \multicolumn{2}{c|}{random init}\\
\hline
img~+~$\sigma$ & KSVD & PGM &  MSM & PGM &  MSM  \\
\hline
walkbridge + 5 & 35.70 & 35.71 & \cthree{35.71} & \ctwo{35.72} & \cone{35.75} \\
walkbridge + 10 & 31.07 & \cthree{31.07} & \ctwo{31.17} & 31.07 & \cone{31.17} \\
walkbridge + 20 & 27.01 & \cthree{27.11} & \ctwo{27.21} & 27.07 & \cone{27.23} \\
walkbridge + 30 & 24.93 & 25.08 & \cone{25.21} & \cthree{25.09} & \ctwo{25.19} \\
walkbridge + 40 & 23.71 & \cthree{23.85} & \ctwo{23.87} & 23.84 & \cone{23.90} \\
mandrill + 5 & 35.18 & \ctwo{35.21} & 35.20 & \cone{35.22} & \cthree{35.21} \\
mandrill + 10 & 30.36 & 30.38 & \ctwo{30.47} & \cthree{30.38} & \cone{30.47} \\
mandrill + 20 & 26.02 & \cthree{26.15} & \ctwo{26.31} & 26.12 & \cone{26.33} \\
mandrill + 30 & 23.75 & \cthree{23.96} & \cone{24.19} & 23.96 & \ctwo{24.18} \\
mandrill + 40 & 22.37 & \cthree{22.61} & \ctwo{22.78} & 22.57 & \cone{22.81} \\
cameraman + 5 & 40.17 & 40.38 & \cone{40.78} & \cthree{40.43} & \ctwo{40.72} \\
cameraman + 10 & 36.04 & 36.07 & \cone{36.52} & \cthree{36.10} & \ctwo{36.50} \\
cameraman + 20 & \cthree{32.03} & 31.77 & \ctwo{32.11} & 31.68 & \cone{32.13} \\
cameraman + 30 & \cone{29.91} & 29.30 & \cthree{29.56} & 29.18 & \ctwo{29.63} \\
cameraman + 40 & \cone{28.39} & 27.55 & \ctwo{27.80} & 27.56 & \cthree{27.78} \\
livingroom + 5 & \cthree{37.00} & 36.97 & \cone{37.10} & 36.94 & \ctwo{37.07} \\
livingroom + 10 & 32.98 & \cthree{33.02} & \ctwo{33.19} & 32.92 & \cone{33.27} \\
livingroom + 20 & 29.22 & 29.16 & \ctwo{29.55} & \cthree{29.23} & \cone{29.57} \\
livingroom + 30 & 27.04 & 27.04 & \ctwo{27.43} & \cthree{27.06} & \cone{27.45} \\
livingroom + 40 & \cthree{25.62} & 25.55 & \ctwo{25.78} & 25.59 & \cone{25.81} \\
\hline
\end{tabular}}\scalebox{\boxscale}{
\begin{tabular}{|p{1.9cm}|p{0.8cm}|p{0.8cm}|p{0.8cm}|p{0.8cm}|p{0.8cm}|}
\hline
 \multicolumn{2}{|c|}{} & \multicolumn{2}{c|}{OMP init} & \multicolumn{2}{c|}{random init}\\
\hline
img~+~$\sigma$ & KSVD & PGM &  MSM & PGM &  MSM  \\
\hline
lake + 5 & \cone{36.77} & 36.74 & \ctwo{36.77} & 36.73 & \cthree{36.75} \\
lake + 10 & \cthree{32.84} & 32.75 & \ctwo{32.86} & 32.77 & \cone{32.90} \\
lake + 20 & \cthree{29.32} & 29.19 & \ctwo{29.36} & 29.23 & \cone{29.38} \\
lake + 30 & \ctwo{27.32} & 27.04 & \cthree{27.30} & 27.10 & \cone{27.33} \\
lake + 40 & \cone{25.94} & 25.66 & \cthree{25.80} & 25.59 & \ctwo{25.80} \\
lena + 5 & \ctwo{38.30} & 38.22 & \cone{38.31} & 38.22 & \cthree{38.29} \\
lena + 10 & 34.81 & 34.79 & \cone{34.97} & \cthree{34.83} & \ctwo{34.91} \\
lena + 20 & \cone{31.48} & 31.16 & \cthree{31.30} & 31.18 & \ctwo{31.34} \\
lena + 30 & \cone{29.50} & 28.96 & \ctwo{29.24} & 28.91 & \cthree{29.11} \\
lena + 40 & \cone{28.07} & 27.35 & \cthree{27.52} & 27.43 & \ctwo{27.57} \\
blonde + 5 & 36.98 & \cthree{37.00} & \ctwo{37.06} & 37.00 & \cone{37.08} \\
blonde + 10 & 33.23 & \cthree{33.31} & \ctwo{33.37} & 33.27 & \cone{33.43} \\
blonde + 20 & \cthree{29.83} & 29.78 & \ctwo{29.99} & 29.75 & \cone{30.00} \\
blonde + 30 & \ctwo{28.00} & 27.77 & \cone{28.04} & 27.73 & \cthree{27.95} \\
blonde + 40 & \cone{26.82} & 26.33 & \cthree{26.41} & 26.31 & \ctwo{26.63} \\
barbara + 5 & \cthree{37.50} & 37.46 & \cone{37.74} & 37.42 & \ctwo{37.71} \\
barbara + 10 & 33.36 & \cthree{33.37} & \ctwo{33.65} & 33.33 & \cone{33.67} \\
barbara + 20 & 29.19 & 29.22 & \cone{29.70} & \cthree{29.30} & \ctwo{29.63} \\
barbara + 30 & 26.79 & 26.84 & \ctwo{27.30} & \cthree{26.92} & \cone{27.36} \\
barbara + 40 & 25.11 & 25.02 & \ctwo{25.71} & \cthree{25.16} & \cone{25.83} \\
\hline
\end{tabular}}\scalebox{\boxscale}{
\begin{tabular}{|p{1.9cm}|p{0.8cm}|p{0.8cm}|p{0.8cm}|p{0.8cm}|p{0.8cm}|}
\hline
 \multicolumn{2}{|c|}{} & \multicolumn{2}{c|}{OMP init} & \multicolumn{2}{c|}{random init}\\
\hline
img~+~$\sigma$ & KSVD &  PGM & MSM  &  PGM & MSM \\
\hline
boat + 5 & \cthree{36.94} & 36.94 & \ctwo{36.98} & 36.93 & \cone{37.01} \\
boat + 10 & \cthree{33.10} & 33.02 & \ctwo{33.24} & 33.05 & \cone{33.31} \\
boat + 20 & 29.47 & 29.40 & \cone{29.66} & \cthree{29.52} & \ctwo{29.65} \\
boat + 30 & \ctwo{27.50} & 27.27 & \cthree{27.50} & 27.29 & \cone{27.52} \\
boat + 40 & \cone{26.16} & 25.84 & \ctwo{26.13} & 25.86 & \cthree{26.03} \\
pirate + 5 & \cthree{36.49} & 36.43 & \cone{36.54} & 36.42 & \ctwo{36.50} \\
pirate + 10 & \cthree{32.19} & 32.09 & \ctwo{32.25} & 32.10 & \cone{32.29} \\
pirate + 20 & \cthree{28.33} & 28.23 & \ctwo{28.42} & 28.24 & \cone{28.44} \\
pirate + 30 & \ctwo{26.31} & 26.12 & \cthree{26.28} & 26.13 & \cone{26.33} \\
pirate + 40 & \cone{24.94} & 24.71 & \cthree{24.86} & 24.72 & \ctwo{24.86} \\
house + 5 & 38.78 & 38.76 & \cone{38.91} & \cthree{38.81} & \ctwo{38.87} \\
house + 10 & 34.99 & 34.97 & \ctwo{35.07} & \cthree{35.05} & \cone{35.10} \\
house + 20 & \cone{31.83} & 31.36 & \ctwo{31.53} & 31.34 & \cthree{31.53} \\
house + 30 & \cone{29.79} & 29.22 & \cthree{29.27} & 29.04 & \ctwo{29.31} \\
house + 40 & \cone{28.16} & 27.40 & \ctwo{27.63} & 27.37 & \cthree{27.46} \\
jetplane + 5 & 38.84 & 38.87 & \cone{39.06} & \cthree{38.93} & \ctwo{39.05} \\
jetplane + 10 & 34.98 & 34.99 & \ctwo{35.17} & \cthree{35.01} & \cone{35.22} \\
jetplane + 20 & \cone{31.30} & 31.04 & \cthree{31.28} & 31.01 & \ctwo{31.30} \\
jetplane + 30 & \cone{29.11} & 28.64 & \cthree{28.87} & 28.68 & \ctwo{28.95} \\
jetplane + 40 & \cone{27.57} & 27.05 & \cthree{27.19} & 26.97 & \ctwo{27.26} \\
\hline
\end{tabular}
}
\caption{Comparisons of SNR values for the sparse coding based image denoising problem with OMP initialization and random initialization. The $1^{st}$, $2^{nd}$, and $3^{rd}$ best results are colored with \cone{red}, \ctwo{blue} and \cthree{green}, respectively. }
\label{table:snr:obj}
\end{table*}

\subsection{Cardinality Regularized Sparse Coding}
Sparse coding is a popular unsupervised feature learning technique for data representation that is widely used in computer vision and medical imaging. Motivated by recent success in $\ell_0$ norm modeling \cite{yuan2015l0tv,BaoJQS16,Yang2016}, we consider the following cardinality regularized (\ie $\ell_0$ norm) sparse coding problem:
\beq \label{eq:card:sparse:coding}
\textstyle \underset{\bbb{W},\bbb{H}}{\min}~~\frac{1}{2}\|\bbb{Y}-\bbb{WH}\|_F^2 + \lambda \|\bbb{H}\|_0,~s.t.~\|\bbb{W}(:,i)\|_{2}=1 ~\forall i
\eeq
\noi with $\bbb{W}\in \mathbb{R}^{m \times n}$ and $\bbb{H}\in \mathbb{R}^{n \times d}$. Existing solutions for this problem are mostly based on the family of proximal point methods \cite{nesterov2013introductory,nesterov2013introductory,BaoJQS16}. We compare MSM with the following methods: (1) Proximal Gradient Method (PGM) with constant step size, (2) PGM with line search, (3) accelerated PGM with constant step size, and (4) accelerated PGM with line search.

We evaluate all the methods for the application of image denoising. Following \cite{aharon2006img,BaoJQS16}, we set the dimension of the dictionary to $n = 256$. The dictionary is learned from $m=1000$ image patches randomly chosen from the noisy input image. The patch size is $8\times 8$, leading to $d=64$. The experiments are conducted on 16 conventional test images with different noise standard deviations $\sigma$. For the regulization parameter $\lambda$, we sweep over $\{1,3,5,7,9,...,10000,30000,50000,70000,90000$$\}$.

First, we compare the objective values for all methods by \emph{fixing} the variable $\bbb{W}$ to an over-complete DCT dictionary \cite{aharon2006img} and \emph{only} optimizing over $\bbb{H}$. We compare all methods with varying regularization parameter $\lambda$ and different initial points that are either generated by random Gaussian sampling or the Orthogonal Matching Pursuit (OMP) method \cite{tropp2007signal}. In Figure \ref{fig:convergece:obj}, we observe that MSM converges rapidly in 10 iterations. Moreover, it often generates much better local optimal solutions than the compared methods.

Second, we evaluate the methods according to Signal-to-Noise Ratio (SNR) value wrt the groundtruth denoised image. In this case, we minimize over $\bbb{W}$ and $\bbb{H}$ \emph{alternatingly} with different initial points (OMP initialization or standard normal random initialization). For updating $\bbb{W}$, we use the same proximal gradient method as in \cite{BaoJQS16}. For updating $\bbb{H}$, since the accelerated PGM does not necessarily present better performance than canonical PGM and since the line search strategy does not present better performance than the simple constant step size, we only compare with PGM, which has been implemented in \cite{BaoJQS16} \footnote{Code: \url{http://www.math.nus.edu.sg/~matjh/research/research.htm}} and KSVD \cite{aharon2006img} \footnote{Code: \url{http://www.cs.technion.ac.il/~elad/software/}}. In our experiments, we find that MSM achieves lower objectives than PGM in all cases. We do not report the objective values here but only the best SNR value, since (i) the best SNR result does not correspond to the same $\lambda$ for PGM and MSM, and (ii) KSVD does not solve exactly the same problem in (\ref{eq:card:sparse:coding})\footnote{In fact, it solves an $\ell_0$ norm constrained problem using a greedy pursuit algorithm and performs a codebook update using SVD. It may not necessarily converge, which motivates the use of the alternating minimization algorithm in \cite{BaoJQS16}. }. We observe that MSM is generally 4-8 times faster than KSVD. This is not surprising, since KSVD needs to call OMP to update the dictionary $\bbb{H}$, which involves high computational complexity while our method only needs to call a generalized Gaussian elimination procedure in each iteration. In Table \ref{table:snr:obj}, we summarize the results, from which we make two conclusions. (i) The two initialization strategies generally lead to similar SNR results. (ii) Our MSM method generally leads to a larger SNR than PGM and a comparable SNR as KSVD, but in less time.

%% file: conclusion.tex
\section{Conclusions and Future Work}

This paper presents a matrix splitting method for composite function minimization. We rigorously analyze its convergence behavior both in convex and non-convex settings. Experimental results on nonnegative matrix factorization and cardinality regularized sparse coding demonstrate that our methods achieve state-of-the-art performance.

Our future work focuses on several directions. (i) We will investigate the possibility of further accelerating the matrix splitting method by Nesterov's momentum strategy \cite{nesterov2013introductory} or Richardson's extrapolation strategy. (ii) It is interesting to extend the classical sparse Gaussian elimination \cite{saad2003iterative} technique to compute the triangle proximal operator in our method, which is expected to be more efficient when the matrix $\bbb{A}$ is sparse. (iii) We are interested in incorporating the proposed algorithms into alternating direction method of multipliers \cite{He2012,boyd2011distributed} as an alternative solution to the existing proximal/linearized method.